%% file: LapFiltersTAP.tex
\documentclass[journal,xcolor=dvipsnames]{IEEEtran}
\ifCLASSINFOpdf
\else
\fi
\input{Packages.tex}

\input{Macros.tex}

\input{preamble}

\usepackage[bb=dsserif]{mathalpha}

\input pdf-trans
\newbox\qbox
\def\usecolor#1{\csname\string\color@#1\endcsname\space}
\newcommand\bordercolor[1]{\colsplit{1}{#1}}
\newcommand\fillcolor[1]{\colsplit{0}{#1}}
\newcommand\outline[1]{\leavevmode%
  \def\maltext{#1}%
  \setbox\qbox=\hbox{\maltext}%
  \boxgs{Q q 2 Tr \thickness\space w \fillcol\space \bordercol\space}{}%
  \copy\qbox%
}
\makeatother
\newcommand\colsplit[2]{\colorlet{tmpcolor}{#2}\edef\tmp{\usecolor{tmpcolor}}%
  \def\tmpB{}\expandafter\colsplithelp\tmp\relax%
  \ifnum0=#1\relax\edef\fillcol{\tmpB}\else\edef\bordercol{\tmpC}\fi}
\def\colsplithelp#1#2 #3\relax{%
  \edef\tmpB{\tmpB#1#2 }%
  \ifnum `#1>`9\relax\def\tmpC{#3}\else\colsplithelp#3\relax\fi
}
\bordercolor{black}
\fillcolor{white}
\def\thickness{.3}
\newcommand{\dSigma}{\text{\outline{$\mat \Sigma$}}}
\newcommand{\dLambda}{\text{\outline{$\mat \Lambda$}}}
\renewcommand{\matd}[1]{\text{\outline{$\mat #1$}}}

\begin{document}
\title{Laplacian Filtered Loop-Star Decompositions and Quasi-Helmholtz Laplacian Filters:\\Definitions, Analysis, and Efficient Algorithms}

\author{Adrien~Merlini,~\IEEEmembership{Member,~IEEE,}
        Clément~Henry,~\IEEEmembership{Member,~IEEE,}
        Davide~Consoli,~\IEEEmembership{Student Member,~IEEE,}
        Lyes~Rahmouni,
        Alexandre~Dély,
        and~Francesco~P.~Andriulli,~\IEEEmembership{Senior Member,~IEEE}% <-this % stops a space
\thanks{D. Consoli, L. Rahmouni, A. Dély, and F. P. Andriulli are with the Department of Electronics and Telecommunications, Politecnico di Torino, 10129 Torino, Italy; e-mail: name.surname@polito.it.}% <-this % stops a space
\thanks{A. Merlini and C. Henry are with the Microwave department, IMT Atlantique, 29238 Brest cedex 03, France; e-mail: name.surname@imt-atlantique.fr.}% <-this % stops a space
\thanks{Manuscript received April 19, 2005.}}


\markboth{Journal of \LaTeX\ Class Files,~Vol.~14, No.~8, August~2015}%
{Shell \MakeLowercase{\textit{et al.}}: Bare Demo of IEEEtran.cls for IEEE Journals}



\maketitle

\begin{abstract}
Quasi-Helmholtz decompositions are fundamental tools in integral equation modeling of electromagnetic problems because of their ability of rescaling solenoidal and non-solenoidal components of solutions, operator matrices, and radiated fields. These tools are however incapable, \textit{per se}, of modifying the refinement-dependent spectral behavior of the different operators and often need to be combined with other preconditioning strategies. 
This paper introduces the new concept of  filtered quasi-Helmholtz decompositions proposing them in two incarnations: the filtered Loop-Star functions and the quasi-Helmholtz Laplacian filters. Because they are capable of manipulating large parts of the operators' spectra, new families of preconditioners and fast solvers can be derived from these new tools. A first application to the case of the frequency and $h$-refinement preconditioning of the electric field integral equation
is presented together with numerical results showing the practical effectiveness of the  newly  proposed decompositions.
\end{abstract}


% \fi
\IEEEpeerreviewmaketitle

\section{Introduction}

\IEEEPARstart{I}{ntegral} equation formulations are effective numerical strategies for modeling radiation and scattering by perfectly electrically conducting objects \cite{nedelec_acoustic_2001,gibson_method_2014,jin_theory_2015}. Their effectiveness primarily derives from the fact that they only require the scatterers' surfaces to be discretized, automatically impose radiation conditions and, thanks to the advent of fast algorithms \cite{chew_fast_2001}, give rise to linear-in-complexity approaches when solved with iterative schemes---provided that the conditioning of the linear system matrices resulting from their discretizations is independent of the number of unknowns \cite{axelsson_iterative_1996}. Among the well-established formulations, the electric field integral equation (EFIE) plays a crucial role, both in itself and within combined field formulations \cite{colton_integral_2013}. The EFIE, lamentably, becomes ill-conditioned when the frequency is low or the discretization density high \cite{adrianElectromagneticIntegralEquations2021}. These phenomena---respectively known as the low-frequency and $h$-refinement breakdowns---cause the solution of the EFIE to become increasingly challenging to obtain, as the number of iterations of the solution process grows unbounded, which jeopardizes the possibility of achieving an overall linear complexity. 

Traditional approaches to tackle the low-frequency breakdown rely on standard quasi-Helmholtz decompositions such as Loop-Star/Tree bases \cite{wilton_topological_1983,vecchi_loop-star_1999,lee_loop_2003,eibert_iterative-solver_2004,andriulli_loop-star_2012} that, despite curing the low-frequency behavior, worsen the $h$-refinement ill-conditioning of the EFIE \cite{eibert_iterative-solver_2004} because of the derivative nature of the change of basis \cite{andriulli_loop-star_2012}.
A way to circumvent the issue is the use of hierarchical strategies both on structured  \cite{vipianaMultiresolutionMethodMoments2005,andriulli_multiresolution_2007} and unstructured meshes \cite{andriulli_hierarchical_2008,chenMultiresolutionCurvilinearRao2009,guzman_hierarchical_2016-1,adrian_hierarchical_2017}. These schemes, when designed properly, can solve both the low-frequency and the $h$-refinement problems but still rely on the construction on an explicit, basis-based, quasi-Helmholtz decomposition that requires the cumbersome detection of topological loops whenever handles are present in the geometry. A popular alternative strategy leverages Calderón identities to form a second kind integral equation out of the EFIE. Calderón approaches concurrently solve the low-frequency and the $h$-refinement breakdowns without calling for an explicit quasi-Helmholtz decomposition \cite{adams_stabilisation_1999,contopanagos_well-conditioned_2002,adams_physical_2004,adams_numerical_2004,borel_new_2005,andriulli_multiplicative_2008,stephanson_preconditioned_2009,borel2005new}. 
In their standard incarnations they do, however, require the use of a dual discretization and global loop handling, because global loops reside in the static null-space of the Calderón operator. The introduction of implicit quasi-Helmholtz decompositions via the so called quasi-Helmholtz projectors \cite{andriulli_well-conditioned_2013}, when combined with Calderón approaches, led to the design of several well-conditioned formulations, free from static nullspaces (see \cite{andriulli_well-conditioned_2013,dobbelaere_calderon_2015, merlini_magnetic_2020,adrianElectromagneticIntegralEquations2021} and references therein) and, in some incarnations, free from the need of performing a barycentric refinement \cite{adrian_refinement-free_2018}. Quasi-Helmholtz projectors have shown to be an effective and efficiently computable tool for performing quasi-Helmholtz decompositions, but, by themselves, they can only tackle the low-frequency breakdown and must be combined with Calderón-like strategies that involve multiple operators, to obtain $h$-refinement spectral preconditioning effects.
A set of tools as versatile as the projectors that could also manipulate the operator spectra beyond a simple rescaling would thus be desirable.

This paper introduces  such a new family of tools. The contribution of this work is in fact threefold: (i) we will introduce the concept of Laplacian-filtered Loop-Star decompositions, a new quasi-Helmholtz decomposition approach that will allow for a finer tuning of the operator spectrum with respect to their standard Loop-Star counterparts. 
(ii) Just like standard Loop-Star bases give rise to the quasi-Helmholtz projectors, a suitable choice of projections on the filtered Loop-Star spaces will give rise to a new family of mathematical objects, the quasi-Helmholtz Laplacian filters, that will not require handling basis functions and explicit decompositions, while still providing the spectral tuning properties of (i). (iii) We will obtain new frequency and $h$-refinement preconditioners for the EFIE based on (i) and (ii) that represent a natural first application of the newly proposed techniques. We believe, however, that the applicability of the new spectral filters will extend beyond EFIE preconditioning in further investigations.
The contribution will be further enriched by a section devoted to efficient implementations of the newly defined tools that will be obtained by leveraging strategies developed in the context of polynomial preconditioning approaches and in graph wavelet theory \cite{johnsonPolynomialPreconditionersConjugate1983,ashbyComparisonAdaptiveChebyshev1992,hammond_wavelets_2011,levie2018cayleynets}. Numerical results will then corroborate and confirm our theoretical considerations.

The paper is organized as follows: the background material and the notation are presented in Section~\ref{sec:background}, the new Laplacian filtered Loop-Star decompositions and the quasi-Helmholtz Laplacian filters are presented in Section~\ref{sec:filtered_LS_bases} and Section~\ref{sec:projector_filters}, respectively, along with their main properties. Various strategies for computing the filters, in practical scenarios, are detailed in Section~\ref{sec:filtering_algos}.
Preconditioners tackling simultaneously the low-frequency and $h$-refinement breakdowns of the EFIE are then derived in Section~\ref{sec:EFIE_application}. Finer details relating to the implementation and computation of the filters and preconditioners are then presented in Section~\ref{eq:implementation_details}. Illustrations of the effectiveness of the schemes are provided in Section~\ref{sec:numerical_results}, before concluding in Section~\ref{sec:conclusion}. 

Preliminary results from this work were presented in the conference contributions \cite{rahmouni_new_2019,9886850}.

%\todo[color=green,inline]{FPA: cite gribonval and polinomial preconditioning}
%\todo[inline,color=cyan]{FPA: This comments applies also to the introduction. We should find here a polynomial preconditioning approach that predates the wavelet guys already using such an approach.}
%\todo[inline]{AM: change the notation for $\mat T_\mathrm{h}$, $\tmat \Sigma_p$, and mesh $h$.}


%\td{Final and very definitive plan of the paper:
%Define the filtered loop and star matrices\newline
%Define the filtered Graph Laplacian projectors\newline
%Define the algorithms to achieve them\newline
%Test and verfy the laplacian filters schemes in a matlab environment\newline
%Define the projectors with gram matrices\newline
%Show that they do diagonalize the vector laplace beltrami operator, once used as wavelets\newline
%Test the surrogate of the EFIE\newline
%Verify on paper that the preconditioner diagonalize the surrogate\newline
%Show that they can be used to precondition the EFIE
%}
\section{Notation and Background}\label{sec:background}

\begin{figure}
        \centering
        \includegraphics{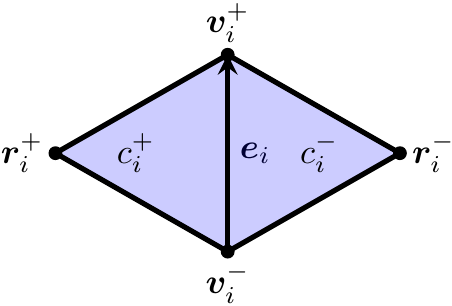}
        \caption{Convention used for the RWGs: each function is defined on the two triangles $c_i^+$ and $c_i^-$ that are formed with their common edge $\vt e_i$ and the vertices $\vr_i^+$ and $\vr_i^-$, respectively.}
        \label{fig:RWGconv}
\end{figure}

Let $\Gamma$ be a smooth surface modeling the boundary of a perfectly electrically conducting (PEC), closed scatterer enclosed in a homogeneous background medium with permittivity $\epsilon$ and permeability $\mu$. The boundary $\Gamma$ can be multiply connected and contain holes. We denote by $\uv{n}(\vr)$ the outward pointing normal field at $\vr$. When illuminated by a time-harmonic incident electric field $\vt{E}^i$, a surface current density $\vt{J}$ is induced on $\Gamma$ that satisfies the electric field integral equation (EFIE)
\begin{equation}
    \label{eq:EFIE}
    \vecop{T} \vt{J}=\vecop{T}_\mathrm{s} \vt{J} +  \vecop{T}_\mathrm{h} \vt{J}=-\uv{n} \times \vt{E^i}\,,
\end{equation}
where
\begin{gather}
    \vecop T_\mathrm{s} \vt{J} = \uv{n}(\vt{r}) \times  \im k  \int_\Gamma \frac{\e^{\im k\| \vt{r} - \vt{r'} \|}}{4\pi \|\vt{r} - \vt{r'} \|}   \vt{J}(\vt{r'})  \dd S(\vt{r'})\,, \\
    \vecop T_\mathrm{h} \vt{J} = -\uv{n}(\vt{r}) \times \frac{1}{\im k} \nabla
    \int_\Gamma \frac{\e^{\im k\| \vt{r} - \vt{r'} \|}}{4\pi\|\vt{r} - \vt{r'} \|}  \nabla' \cdot \vt{J}(\vt{r'})  \dd S(\vt{r'})\,,
\end{gather}
and $k$ is the wavenumber of the electromagnetic wave in the background medium.
Equation~\eqref{eq:EFIE} can be solved numerically by approximating $\Gamma$ with triangular elements of average edge length $h$ and by approximating the current density as $\vt{J} \approx \sum_{n=1}^{N} [\vec j]_n \vt f_n$ with the Rao-Wilton-Glisson basis functions $\{\vt f_n\}_n$ \cite{rao_electromagnetic_1982}, in which $N$ is the number of edges in the mesh, $\vec j$ is the vector of the coefficients of the expansion, and $\vt f_n$ is defined as
\begin{equation} \label{eq:rwg_def}
        \vt{f_n}(\vr) =
        \begin{dcases}
                \phantom{-}\frac{\vr - \vr_n^+}{2 A_n^+} & \text{if } \vr \in c_n^+ \\
                -\frac{\vr - \vr_n^-}{2 A_n^-} & \text{if } \vr \in c_n^-,
        \end{dcases}
\end{equation}
where the notation of Fig.~\ref{fig:RWGconv} was employed and where $A_n^\pm$ is the area of the cell $c_n^\pm$.

The final step to obtain the discretized EFIE is to test \eqref{eq:EFIE} with the rotated RWG functions $\{\uv{n} \times \vt{f_n}\}$, which results in the linear system
\begin{equation} \label{eq:disc_EFIE}
    \mat{T} \vec j = \left( \mat T_\mathrm{s}  + \mat T_\mathrm{h} \right) \vec{j} = \vec{v}\,,
\end{equation}
in which $[\mat T_\mathrm{s} ]_{mn} = \langle \uv{n} \times \vt{f_m} ,  \vecop{T}_s (\vt{f}_n)\rangle$, $[\mat T_\mathrm{h}]_{mn} = \langle \uv{n} \times \vt{f_m} , \vecop{T}_h (\vt{f}_n)\rangle$, $[\vec{v}]_{m} = \langle \uv{n} \times \vt{f_m} , -\uv{n} \times \vt{E^i} \rangle$, and $\langle \vt{a}, \vt{b} \rangle=\int_\Gamma\vt{a}\cdot \vt{b} \dd s$.
The EFIE can also be discretized on the dual mesh using dual functions defined on the barycentric refinement. Both Buffa-Christiansen \cite{buffa_dual_2007} and Chen-Wilton \cite{chen_electromagnetic_1990} elements can be used for this dual discretization. For the sake of brevity, we will omit the explicit definitions of the dual elements that will be denoted by $\{\vt g_n\}_n$ in the following; the reader can refer to \cite{adrianElectromagneticIntegralEquations2021} and references therein for a more detailed treatment. We will also need the definition of the standard and dual Gram matrices whose entries are $\left[\mat{G}\right]_{mn}=\langle \vv{f}_m,\vv{f}_n\rangle$ and $\left[\matd{G}\right]_{mn}=\langle \vv{g}_m,\vv{g}_n\rangle$. While they are not required for the discretization of the EFIE itself, we introduce the patch and pyramids scalar basis functions sets $\{p_n\}$ and $\{\lambda_n\}$, respectively composed of $N_S$ and $N_L$ functions, that will be required for some of the following developments. These functions are defined as
\begin{equation}
    p_m(\vv r) = \begin{dcases}
      A_m^{-1} & \text{if } \vv r \in c_m\,,\\
      0 & \text{otherwise,}
    \end{dcases}
\end{equation}
and
\begin{equation}
     \lambda_m(\vv r) = \begin{dcases}
      1 & \vv r = \vv v_m\,,\\
      0 & \vv r = \vv v_n\,, n \neq m\,,\\
      \text{linear} & \text{otherwise,}
    \end{dcases}
\end{equation}
where $\{\vv v_n\}_n$ are the vertices of the mesh. The number of these functions can be deduced from the mesh properties: $N_S$ is the number of mesh triangles and $N_L$ is the number of mesh vertices. The Gram matrices corresponding to these bases are $\mat G_p$ for the patch functions and $\mat G_\lambda$ for the pyramids with $[\mat G_p]_{mn} = \left<p_m, p_n\right>$ and $[\mat G_\lambda]_{mn} = \left<\lambda_m, \lambda_n\right>$. The dual of these functions, living in the barycentric refinement of the original mesh will also be required, and their definitions, omitted here for conciseness, can be found in \cite{adrian_refinement-free_2018}. The $N_L$ dual patches will be designated as $\{\tilde p_n\}_n$, the $N_S$ dual pyramids as $\{\tilde \lambda_n\}_n$, and the corresponding gram matrices as $\matd G_{\tilde p}$ and $\matd G_{\tilde \lambda}$ with $[\matd G_{\tilde p}]_{mn} = \left<\tilde p_m, \tilde p_n\right>$ and $[\matd G_{\tilde\lambda}]_{mn} = \left<\tilde \lambda_m, \tilde \lambda_n\right>$.

Because this contribution deals with discrete quasi-Helmholtz decompositions we will recall some of their properties.
The continuous solution $\vv{J}$ can be decomposed into a solenoidal, irrotational, and (on non simply-connected manifolds) harmonic components
\begin{equation}\label{eq:Hdecomposition}
\vv{J}=\nabla \times \uv{n}\lambda + \nabla_s \phi + \vv{h}\,.
\end{equation}
When $\vv{J}$ is discretized by the approximate expansion in RWG functions, a discrete counterpart of \eqref{eq:Hdecomposition} holds for the coefficient vector $\vec{j}$
\begin{equation}\label{eq:discreteHdecomposition}
    \vec{j}= \mat{\Lambda}\vec{l}+\mat{\Sigma}\vec{\sigma}+\mat{H}\vec{h}\,,
\end{equation}
where $\mat{\Lambda} \in \mathbb{R}^{N\times N_L}$ and $\mat{\Sigma} \in \mathbb{R}^{N\times N_S}$ are the Loop-to-RWG and Star-to-RWG transformation matrices \cite{mautzEfieldSolutionConducting1984,limNovelTechniqueCalculate,wuStudyTwoNumerical1995,vecchi_loop-star_1999,andriulli_loop-star_2012} defined, following the convention of Fig.~\ref{fig:RWGconv} and the definition of the RWG in \eqref{eq:rwg_def}, as
\begin{equation}\label{eq:loopmat}
  \left[\mat{\Lambda}\right]_{mn}=\begin{dcases}
    \phantom{-}1& \text{if  node } n \text{ equals } v_m^+\,,\\
    -1& \text{if  node } n \text{ equals } v_m^-\,,\\
    \phantom{-}0& \text{otherwise,}
  \end{dcases}
\end{equation}
and
\begin{equation}\label{eq:starmat}
  \left[\mat{\Sigma}\right]_{mn}=\begin{dcases}
    \phantom{-}1& \text{if the cell } n \text{ equals } c_m^+\,,\\
    -1& \text{if the cell } n \text{ equals } c_m^-\,,\\
    \phantom{-}0& \text{otherwise.}
  \end{dcases}
\end{equation}
With these definitions $\mat\Lambda^\T\mat\Lambda$ and $\mat\Sigma^\T\mat\Sigma$ are respectively the vertices- and the cells-based graph Laplacians \cite{andriulli_loop-star_2012}.
The explicit use of the change of basis matrix $\mat{H}$ will not be required and we omit here its explicit definition, for the sake of conciseness, which could however be found in \cite{adrianElectromagneticIntegralEquations2021} and references therein.

With the definitions above, the standard quasi-Helmholtz projectors
\cite{andriulli_loop-star_2012,andriulli_well-conditioned_2013} are defined as
\begin{align}
    \begin{split}
    \mat P^\Sigma &=\mat\Sigma \left(\mat\Sigma^\T\mat\Sigma\right)^{+} \mat\Sigma^\T\label{eq:sigmaprojdef}\,,\\
    \mat P^{\Lambda H} &=\mat I-  \mat P^\Sigma
    \end{split}
\end{align}
for the primal ones,
%\begin{align}
%       \matd{P}^\Lambda &=\mat\Lambda \left(\mat\Lambda^\T\mat\Lambda\right)^{+} \mat\Lambda^\T\label{eq:lambdaprojdef}\,,\\
%  \matd{P}^{\Sigma H} &=\mat I-  \matd{P}^\Lambda\,.
%\end{align}
\begin{align}
    \begin{split}
       \matd{P}^\Lambda &=\mat\Lambda \left(\mat\Lambda^\T\mat\Lambda\right)^{+} \mat\Lambda^\T\label{eq:lambdaprojdef}\,,\\
        \matd{P}^{\Sigma H} &=\mat I-  \matd{P}^\Lambda
    \end{split}
\end{align}
for the dual ones, and
\begin{equation}
    \mat P^{H} \left(=\matd P^H\right)= \mat I -  \mat P^\Sigma - \matd P^\Lambda\\
\end{equation}
for the projector to quasi-harmonic subspace, where $^+$ denotes the Moore-Penrose pseudo-inverse.

\section{Laplacian Filtered Loop-Star Decompositions} \label{sec:filtered_LS_bases}

In this section, we will extend the notion of Loop-Star bases by introducing the concept of filtered (generating) functions. We will first treat graph-based decompositions (a direct generalization of the standard case) and we will then move on to their Gram matrix normalized counterparts that will be more effective in treating problems involving inhomogeneous meshes.

%\subsection{$l^2$-Orthogonal Decomposition}
\subsection{The Standard Case}

Consider the singular value decomposition (SVD) \cite{golub_matrix_2012} of a matrix $\mat X \in \mathbb{R}^{N\times N_{x}}$
\begin{equation} \label{eq:svd}
    \mat X=\mat U_{\mat X} \mat S_{\mat X} \mat V_{\mat X}^\T
\end{equation}
where $\mat X$ is a placeholder for either $\mat \Sigma$ or $\mat \Lambda$, $\mat U_{\mat X} \in \mathbb{R}^{N\times N}$, $\mat V_{\mat X} \in \mathbb{R}^{N_{x}\times N_{x}}$, and $\mat S_{\mat X} \in \mathbb{R}^{N\times N_{x}}$. The matrices $\mat U_{\mat X}$ and $\mat V_{\mat X}$ are unitary and $\mat S_{\mat X}$ is a block diagonal matrix with the singular values $\sigma_{\mat X,i}$ of $\mat X$ as entries (in decreasing order). Clearly the SVD of $(\mat X^\T\mat X)$ is $\mat V_{\mat X} \mat S_{\mat X}^\T \mat S_{\mat X} \mat V_{\mat X}^\T$, and, by defining the diagonal matrix $\mat L_{\mat X,n} \in \mathbb{R}^{N_{x}\times N_{x}}$, with $1 \leq n \leq N_{x}$, such that
%\td{Be careful because for Loops on open structures we many have the above index to start running from 0 instead of 1.}
%
\begin{equation} \label{eq:def_L}
  \left[\mat L_{\mat X,n} \right]_{ii} =
  \begin{dcases}
    \sigma_{\mat X,i} & \text{if } i > N_{x}-n\,,\\
   \phantom{-}0 & \text{otherwise,}
  \end{dcases}
\end{equation}%
we define the filtered graph Laplacians
\begin{equation}\label{eq:svdX}
(\mat X^\T\mat X)_{n} \coloneqq \mat V_{\mat X} \mat L_{\mat X,n}^2 \mat V_{\mat X}^\T\,,
\end{equation}
from which we introduce the filtered Loop-to-RWG and filtered Star-to-RWG matrices we propose in this work
\begin{align}
    \mat\Sigma_{n}&=\mat\Sigma\left(\mat\Sigma^\T\mat\Sigma\right)^+\left(\mat\Sigma^\T\mat\Sigma\right)_{n}\,,\label{eq:filtered1}\\
    \mat\Lambda_{n}&=\mat\Lambda\left(\mat\Lambda^\T\mat\Lambda\right)^+\left(\mat\Lambda^\T\mat\Lambda\right)_{n}\,.\label{eq:filtered2}
\end{align}
These matrices contains the coefficients of sets of linearly dependent filtered Loop-Star functions.

\subsubsection*{Properties}\label{subsec:LS_props}
We now study some properties of the filtered Loop-Star matrices. Because $\mat\Sigma^\T\mat\Lambda=\matO$ \cite{adrianElectromagneticIntegralEquations2021}, we have $\forall n,m$
\begin{multline} \label{eq:orthogonality1}
    \mat\Sigma_{n}^\T\mat\Lambda_{m} =\left(\mat\Sigma^\T\mat\Sigma\right)_{n}\left(\mat\Sigma^\T\mat\Sigma\right)^+\mat\Sigma^\T\mat\Lambda\left(\mat\Lambda^\T\mat\Lambda\right)^+\left(\mat\Lambda^\T\mat\Lambda\right)_{m}\\
=\matO\,.
\end{multline}
Otherwise said, the new filtered Loop-Star functions are coefficient-orthogonal ($l^2$-orthogonal) like their non-filtered, standard counterparts.

From the definition of $\mat L_{\mat X,n}$ in \eqref{eq:def_L}, it follows that $\mat L_{\mat X,n}\mat L_{\mat X,m}=\mat L^2_{\mat X,\min \{n,m\}}$. Thus from \eqref{eq:svdX}
\begin{align}\label{eq:passagesmin}
   & (\mat X^\T\mat X)_{n} (\mat X^\T\mat X)_{m} =\mat V_{\mat X} \mat L_{\mat X,n} \mat V_{\mat X}^\T\mat V_{\mat X} \mat L_{\mat X,m} \mat V_{\mat X}^\T\nonumber\\
    &=\mat V_{\mat X} \mat L_{\mat X,n} \mat L_{\mat X,m} \mat V_{\mat X}^\T
    =\mat V_{\mat X} \mat L_{\mat X,\min \{n,m\}} \mat L_{\mat X,\min \{n,m\}} \mat V_{\mat X}^\T\nonumber\\
    &=\mat V_{\mat X} \mat L_{\mat X,\min \{n,m\}} \mat V_{\mat X}^\T\mat V_{\mat X} \mat L_{\mat X,\min \{n,m\}} \mat V_{\mat X}^\T\\
     &=(\mat X^\T\mat X)^2_{\min \{n,m\}}\,. \nonumber
    \end{align}
We thus have
\begin{align}
    \mat\Sigma_{m}^\T\mat\Sigma_{n}&=
      \left(\mat\Sigma^\T\mat\Sigma\right)_{m} \left(\mat\Sigma^\T\mat\Sigma\right)^+   \mat\Sigma^\T
    \mat\Sigma\left(\mat\Sigma^\T\mat\Sigma\right)^+\left(\mat\Sigma^\T\mat\Sigma\right)_{n}\nonumber\\
    &=  \left(\mat\Sigma^\T\mat\Sigma\right)^+\left(\mat\Sigma^\T\mat\Sigma\right)_{m}\left(\mat\Sigma^\T\mat\Sigma\right)_{n}\nonumber\\
    &= \left(\mat\Sigma^\T\mat\Sigma\right)^+\left(\mat\Sigma^\T\mat\Sigma\right)_{\min \{n,m\}}\left(\mat\Sigma^\T\mat\Sigma\right)_{\min \{n,m\}}\\
    &=\mat\Sigma_{\min \{n,m\}}^\T\mat\Sigma_{\min \{n,m\}}\,.\nonumber
    \end{align}
Similarly,
\begin{equation}
        \mat\Lambda_{m}^\T\mat\Lambda_{n}=\mat\Lambda_{\min \{n,m\}}^\T\mat\Lambda_{\min \{n,m\}}.
\end{equation}

Given integers such that $m<n<p<q$, the property
\begin{multline}
 \left(\mat\Sigma_{m}-\mat\Sigma_{n}\right)^\T\left(\mat\Sigma_{p}-\mat\Sigma_{q}\right)=\\
    \mat\Sigma_{m}^\T\left(\mat\Sigma_{p}-\mat\Sigma_{q}\right) - \mat\Sigma_{n}^\T\left(\mat\Sigma_{p}-\mat\Sigma_{q}\right)=\matO\,, \label{eq:orthdiffSigma}
\end{multline}
holds and, similarly,
\begin{equation}\label{eq:orthdiffLambda}
    \left(\mat\Lambda_{m}-\mat\Lambda_{n}\right)^\T\left(\mat\Lambda_{p}-\mat\Lambda_{q}\right)=\matO\,.
\end{equation}
Properties \eqref{eq:orthdiffSigma} and \eqref{eq:orthdiffLambda} show that non-intersecting differences of filtered Star or Loop bases are mutually orthogonal (and thus generate linearly independent spaces), a property that will be useful to build invertible changes of basis, as will be shown in Section~\ref{subsec:precond_filtered_LS}.

%\subsection{$L^2$-Orthogonal Functions}
\subsection{Generalization for Non-homogeneously Meshed Geometries}

When the filtered Loop-Star decompositions are to be used on geometries with non-homogenous discretizations, both the standard discretizations of the EFIE and the graph Laplacian matrices may lead to suboptimal performance and a proper normalization with Gram matrices must be employed. In this context, we define the normalized EFIE electromagnetic operator matrices 
\begin{align}
    \tmat{T} &= \mat G^{-1/2} \mat T \mat G^{-1/2}\,,\\
    \tmat T_\mathrm{s} &= \mat G^{-1/2} \mat T_\mathrm{s}  \mat G^{-1/2}\,,\\
    \tmat T_\mathrm{h} &= \mat G^{-1/2} \mat T_\mathrm{h}  \mat G^{-1/2}\,,
\end{align}
and the normalized Loop and Star matrices
%\td{AM: Here the patches in $\mat G_p$ are defined as having $1/A$ value where $A$ is the area of the defining triangle  }
\begin{align}
    \tilde{\mat \Sigma}&=\mat G^{-1/2}\mat\Sigma\mat G_{p}^{1/2}\,,\label{eq:normLS1}\\
    \tilde{\mat \Lambda}&=\mat G^{1/2}\mat\Lambda\mat G_{\lambda}^{-1/2}\,.\label{eq:normLS2}
\end{align}
Following the same strategy as in \eqref{eq:filtered1} and \eqref{eq:filtered2}, the normalized filtered Loop-Star matrices are consistently defined as
\begin{align}
    \mat\SigmaN_{n}&=\mat\SigmaN\left(\mat\SigmaN^\T\mat\SigmaN\right)^+\left(\mat\SigmaN^\T\mat\SigmaN\right)_{n}\,,\label{FiltStarPrimal}\\
    \mat\LambdaN_{n}&=\mat\LambdaN\left(\mat\LambdaN^\T\mat\LambdaN\right)^+\left(\mat\LambdaN^\T\mat\LambdaN\right)_{n}\,. \label{FiltLoopPrimal}
\end{align}

When dealing with dual Loop-Star decomposition matrices, the normalization is different from that of the primal ones, and the dually-normalized Loop and Star transformation matrices are defined as
\begin{align}
    \tilde{\dSigma}&=\matd G^{1/2}\mat\Sigma\matd{G}_{\tilde \lambda}^{-1/2}\,,\\
    \tilde{\dLambda}&=\matd G^{-1/2}\mat\Lambda \matd G_{\tilde p}^{1/2}\,,
\end{align}
and the associated filtered decomposition matrices as
\begin{align}
    \mat\SigmaNd_{n}&=\mat\SigmaNd\left(\mat\SigmaNd^\T\mat\SigmaNd\right)^+\left(\SigmaNd^\T\SigmaNd\right)_{n}\,,\\
    \mat\LambdaNd_{n}&=\mat\LambdaNd\left(\mat\LambdaNd^\T\mat\LambdaNd\right)^+\left(\LambdaNd^\T\LambdaNd\right)_{n}\,.
\end{align}

\subsubsection*{Properties} \label{subsec:norm_LS_props}
The primal and dual normalized Loop-Star bases keep satisfying the orthogonality properties
\begin{align}\label{eq:relorth1}
    \tilde{\mat \Lambda}^\T\tilde{\mat \Sigma}&=\mat G_{\lambda}^{-1/2} \mat\Lambda^\T \mat G^{1/2} \mat G^{-1/2}\mat\Sigma \mat G_{p}^{1/2} = \matO\,,\\
\label{eq:relorth2}
    \tilde{\dSigma}^\T\tilde{\dLambda} &= \matd G_{\tilde p}^{1/2}\mat\Sigma^\T \matd G^{-1/2} \matd G^{1/2} \mat \Lambda \matd G_{\tilde \lambda}^{-1/2} = \matO\,,
\end{align}
because $\mat\Sigma^\T\mat\Lambda=\matO$. Moreover, because \eqref{eq:passagesmin} holds, we obtain, similarly to Section~\ref{subsec:LS_props}, that
\begin{align}
    \tilde{\mat\Sigma}_{m}^\T  \tilde{\mat\Sigma}_{n}  &=\tilde{\mat\Sigma}_{\min \{n,m\}}^\T\tilde{\mat\Sigma}_{\min \{n,m\}}\label{propertiesfirst}\,,\\
    \tilde{\mat\Lambda}_{m}^\T \tilde{\mat\Lambda}_{n} &=\tilde{\mat\Lambda}_{\min \{n,m\}}^\T\tilde{\mat\Lambda}_{\min \{n,m\}}\,,\\
    \tilde{\matd\Lambda}_{m}^\T\tilde{\matd\Lambda}_{n}&=\tilde{\matd\Lambda}_{\min \{n,m\}}^\T\tilde{\matd\Lambda}_{\min \{n,m\}}\,,\\
    \tilde{\matd\Sigma}_{m}^\T \tilde{\matd\Sigma}_{n} &=\tilde{\matd\Sigma}_{\min \{n,m\}}^\T\tilde{\matd\Sigma}_{\min \{n,m\}}\,.\label{propertieslast}
\end{align}
Using these properties and the same reasoning as previously, the counterparts of the properties of the non-normalized filtered Loop-Star matrices can be obtained. In particular the counterparts of \eqref{eq:orthdiffSigma} and \eqref{eq:orthdiffLambda} can be obtained by replacing each matrix with its normalized (``tilde'') counterpart.

\section{Quasi-Helmholtz Laplacian Filters} \label{sec:projector_filters}

Although explicit quasi-Helmholtz decomposition bases are useful in applications in which a direct access to the Helmholtz components of the current is required, oftentimes, especially when the main target is preconditioning and regularization, implicit Helmholtz decompositions can be more efficient. An implicit Helmholtz decomposition was obtained in \cite{andriulli_well-conditioned_2013}, where the concept of quasi-Helmholtz projector was introduced. Following a similar philosophy, and leveraging the filtered Loop-Star functions introduced above, we can now define quasi-Helmholtz Laplacian filters.

\subsection{The Standard Case}
The idea behind the projectors was to obtain a basis-free quasi-Helmholtz decomposition that would not worsen the conditioning of the original equation. If we follow the definitions \eqref{eq:sigmaprojdef} and \eqref{eq:lambdaprojdef} by replacing the standard Star basis with the new filtered sets, we obtain
\begin{multline}
\mat\Sigma_n \left(\mat\Sigma_n^\T\mat\Sigma_n\right)^{+} \mat\Sigma_n^\T=
\mat\Sigma\left(\mat\Sigma^\T\mat\Sigma\right)^+\left(\mat\Sigma^\T\mat\Sigma\right)_{n}\\
\left(\left(\mat\Sigma^\T\mat\Sigma\right)_{n} \left(\mat\Sigma^\T\mat\Sigma\right)^+\mat\Sigma^\T
\mat\Sigma\left(\mat\Sigma^\T\mat\Sigma\right)^+\left(\mat\Sigma^\T\mat\Sigma\right)_{n}\right)^+\\
\left(\mat\Sigma^\T\mat\Sigma\right)_{n} \left(\mat\Sigma^\T\mat\Sigma\right)^+\mat\Sigma^\T
=\mat\Sigma\left(\mat\Sigma^\T\mat\Sigma\right)_{n}^+
\mat\Sigma^\T\,,
\end{multline}
and, similarly,
\begin{equation}
\mat\Lambda_n\left(\mat\Lambda_n^\T\mat\Lambda_n\right)^{+}\mat\Lambda_n^\T=\mat\Lambda\left(\mat\Lambda^\T\mat\Lambda\right)_{n}^+\mat\Lambda^\T\,.
\end{equation}
This justifies the following definitions of the new primal filters
\begin{align}
    \mat P_n^\Sigma&=\mat\Sigma \left(\mat\Sigma^\T\mat\Sigma\right)^{+}_n \mat\Sigma^\T\,,\\
    \mat P_n^{\Lambda H}&=\mat\Lambda \left(\mat\Lambda^\T\mat\Lambda\right)^{+}_n \mat\Lambda^\T +\mat I-\mat P^\Sigma-\matd P^\Lambda\label{eq:defLamH}
\end{align}
and dual filters
\begin{align}
\matd P_n^\Lambda&=\mat\Lambda \left(\mat\Lambda^\T\mat\Lambda\right)^{+}_n \mat\Lambda^\T\,, \\
    \matd P_n^{\Sigma H}&=\mat\Sigma \left(\mat\Sigma^\T\mat\Sigma\right)^{+}_n \mat\Sigma^\T+\mat I-\matd P^\Lambda-\mat P^\Sigma\,.\label{eq:defSigH}
\end{align}
The reader should note that, in the special case of simply connected geometries
$\mat P_n^{\Lambda H}=\mat\Lambda \left(\mat\Lambda^\T\mat\Lambda\right)^{+}_n \mat\Lambda^\T$ and
$\matd P_n^{\Sigma H}=\mat\Sigma \left(\mat\Sigma^\T\mat\Sigma\right)^{+}_n \mat\Sigma^\T$ since
$\matd P^\Lambda+\mat P^\Sigma=\mat I$. Moreover, by construction,
\begin{align}
    \mat P_{N_S}^\Sigma&=\mat P^\Sigma\,,\\
    \matd P_{N_L}^\Lambda&=\matd P^\Lambda\,,
\end{align}
and thus
\begin{align}
    \mat P_{N_L}^{\Lambda H}&= \matd P^\Lambda +\mat I-\mat P^\Sigma-\matd P^\Lambda=\mat P^{\Lambda H}\,,\\
    \matd P_{N_S}^{\Sigma H}&=\mat P^\Sigma+\mat I-\matd P^\Lambda-\mat P^\Sigma= \matd P^{\Sigma H}\,,
\end{align}
which means that, with these definitions, the quasi-Helmholtz Laplacian filters converge to the standard quasi-Helmholtz projectors when the Laplacian is unfiltered ($n=N_{\mat X}$).

\subsubsection*{Properties} \label{subsec:proj_properties}
From these definitions, a few useful properties of the quasi-Helmholtz Laplacian filters can be derived and will be summarized here.
First, the filters still behave as projectors since
\begin{align}\label{eq:proj1}
\begin{split}
    \mat P_n^\Sigma\mat P_n^\Sigma&=\mat\Sigma \left(\mat\Sigma^\T\mat\Sigma\right)^{+}_n
    \mat\Sigma^\T\mat\Sigma \left(\mat\Sigma^\T\mat\Sigma\right)^{+}_n \mat\Sigma^\T\\
    &=\mat\Sigma \left(\mat\Sigma^\T\mat\Sigma\right)^{+}_n
    \mat\Sigma^\T=\mat P_n^\Sigma\,,
        \end{split}\\
    \begin{split}
        \matd P_n^\Lambda\matd P_n^\Lambda&=\mat\Lambda \left(\mat\Lambda^\T\mat\Lambda\right)^{+}_n \mat\Lambda^\T \mat\Lambda \left(\mat\Lambda^\T\mat\Lambda\right)^{+}_n \mat\Lambda^\T \\
    &=\mat\Lambda \left(\mat\Lambda^\T\mat\Lambda\right)^{+}_n \mat\Lambda^\T =\matd P_n^\Lambda\,,
    \end{split}\label{eq:proj2}
\end{align}
and, similarly,
\begin{align}
    \mat P_n^{\Lambda H}\mat P_n^{\Lambda H}&=\mat P_n^{\Lambda H}\,,\label{eq:proj3}\\
    \matd P_n^{\Sigma H}\matd P_n^{\Sigma H}&=\matd P_n^{\Sigma H}\,.\label{eq:proj4}
\end{align}
%The above properties show that quasi-Helmholtz Laplacian Filters are a set of projectors.
Moreover, $\forall m,n$
\begin{multline}
\mat P_{m}^\Sigma\mat P_{n}^{\Lambda H} =
\mat\Sigma \left(\mat\Sigma^\T\mat\Sigma\right)^{+}_n \mat\Sigma^\T
  \mat\Lambda \left(\mat\Lambda^\T\mat\Lambda\right)^{+}_n \mat\Lambda^\T \\
  +\mat\Sigma \left(\mat\Sigma^\T\mat\Sigma\right)^{+}_n \mat\Sigma^\T\left(\mat I-\mat P^\Sigma-\matd P^\Lambda\right)=\matO\,,
\end{multline}
where the properties $\mat\Sigma^\T\mat\Lambda=\matO$ and $\mat \Sigma^\T\left(\mat I-\mat P^\Sigma-\matd P^\Lambda\right)=\matO$ have been used. A similar property and proof hold for the dual projectors
\begin{equation}
\matd P_{m}^\Lambda\matd P_{n}^{\Sigma H}=\matO\,, \quad\forall m,n\,.
\end{equation}

For integers $m<n<p<q$, we have the following orthogonality property
%\mat\Sigma \left(\mat\Sigma^T\mat\Sigma\right)^{+}_n \mat\Sigma^T
\begin{align}
\begin{split}
    \MoveEqLeft[3] \left(\mat P_{m}^\Sigma-\mat P_{n}^\Sigma\right)\left(\mat P_{p}^\Sigma-\mat P_{q}^\Sigma\right)\\
    ={}& \left(\mat\Sigma \left(\mat\Sigma^\T\mat\Sigma\right)^{+}_m \mat\Sigma^\T-\mat\Sigma \left(\mat\Sigma^\T\mat\Sigma\right)^{+}_n \mat\Sigma^\T\right)\\
    &\left(\mat\Sigma \left(\mat\Sigma^\T\mat\Sigma\right)^{+}_p \mat\Sigma^\T-\mat\Sigma \left(\mat\Sigma^\T\mat\Sigma\right)^{+}_q \mat\Sigma^\T\right)\\
    ={}& \mat\Sigma \left(\mat\Sigma^\T\mat\Sigma\right)^{+}_m \mat\Sigma^\T-\mat\Sigma \left(\mat\Sigma^\T\mat\Sigma\right)^{+}_m \mat\Sigma^\T\\
     &+\mat\Sigma \left(\mat\Sigma^\T\mat\Sigma\right)^{+}_n  \mat\Sigma^\T-\mat\Sigma \left(\mat\Sigma^\T\mat\Sigma\right)^{+}_n \mat\Sigma^\T=\matO\,,
     \end{split}
\end{align}
where \eqref{eq:passagesmin} has been used. In a similar way, one can prove that
\begin{align}
      \left(\matd P_{m}^\Lambda-\matd P_{n}^\Lambda\right)\left(\matd P_{p}^\Lambda-\matd P_{q}^\Lambda\right) &=\matO\,.
\end{align}
Moreover, given that $\mat P_n^{\Lambda H}-\mat P_m^{\Lambda H}=\matd P_n^{\Lambda}-\matd P_m^{\Lambda}$ and $\matd P_n^{\Sigma H}-\matd P_m^{\Sigma H}=\mat P_n^{\Sigma}-\mat P_m^{\Sigma}$ $\forall n,m$---which can be deduced from \eqref{eq:defLamH} and \eqref{eq:defSigH}---the remaining properties
\begin{align}
 \left(\mat P_{m}^{\Lambda H}-\mat P_{n}^{\Lambda H}\right)\left(\mat P_{p}^{\Lambda H}-\mat P_{q}^{\Lambda H}\right)&=\matO\,,\\
  \left(\matd P_{m}^{\Sigma H}-\matd P_{n}^{\Sigma H}\right)\left(\matd P_{p}^{\Sigma H}-\matd P_{q}^{\Sigma H}\right)&=\matO \label{eq:lasteq}
\end{align}
follow.
All the properties listed above, will be useful when building invertible transforms, similarly to their basis-based counterpart \eqref{eq:orthdiffSigma}.

\subsection{Generalization for Non-homogeneously Meshed Geometries}
The definitions of the normalized Loop and Star matrices in \eqref{eq:normLS1} 
and \eqref{eq:normLS2}
suggest the following definition for the associated normalized quasi-Helmholtz projectors
\begin{align}
     \Pn^\Sigma&=\SigmaN\left(\SigmaN^\T\SigmaN\right)^{+}\SigmaN^\T\,,\\
    \Pn^{\Lambda}&=\LambdaN \left(\LambdaN^\T\LambdaN\right)^{+}\LambdaN^\T\,.
\end{align}
Moreover, as is proved in Appendix~\ref{app:complementarity}, the complementarity property
\begin{equation}
     \Pn^\Sigma=\mat I -  \Pn^{\Lambda}
\end{equation}
holds on simply connected geometries; on general geometries and together with definitions \eqref{FiltStarPrimal} and \eqref{FiltLoopPrimal}, this justifies the following definition for the normalized quasi-Helmholtz Laplacian filters
\begin{align}
     \Pn_n^\Sigma&=\SigmaN\left(\SigmaN^\T\SigmaN\right)^{+}_n \SigmaN^\T\,,\\
    \Pn_n^{\Lambda H}&=\LambdaN \left(\LambdaN^\T\LambdaN\right)^{+}_n \LambdaN^\T +\mat I-\tilde{\mat P}^\Sigma-\tilde{\mat P}^\Lambda\,.
\end{align}
By analogy, we can define the normalized dual quasi-Helmholtz projectors as
\begin{align}
 \Pnd^{\Lambda}&=\tilde{\dLambda}\left(\tilde{\dLambda}^\T\tilde{\dLambda}\right)^{+}\tilde{\dLambda}^\T\,,\\
     \Pnd^\Sigma&=\tilde{\dSigma}\left(\tilde{\dSigma}^\T\tilde{\dSigma}\right)^{+}\tilde{\dSigma}^\T\,,
\end{align}
with the property
\begin{equation}
    \Pnd^{\Lambda}=\mat I - \Pnd^\Sigma
\end{equation}
holding on simply connected geometries (see Appendix~\ref{app:complementarity} for the proof). Thus, dually to the primal case, we define, on general geometries,
\begin{align}
     \Pnd_n^\Lambda&=\tilde{\dLambda} \left(\tilde{\dLambda}^\T\tilde{\dLambda}\right)^{+}_n\tilde{\dLambda}^\T\,,\\
    \Pnd_n^{\Sigma H}&= \tilde{\dSigma}\left(\tilde{\dSigma}^\T\tilde{\dSigma}\right)^{+}_n \tilde{\dSigma}^\T+\mat I-\Pnd^\Lambda-\Pnd^\Sigma\,.
\end{align}

\subsubsection*{Properties}

Since the primal and dual normalized Loop-Star bases still satisfy the orthogonality properties \eqref{eq:relorth1} and
\eqref{eq:relorth2} and because of the properties \eqref{propertiesfirst}-\eqref{propertieslast},
the same reasoning yields all counterparts of the properties
\eqref{eq:proj1}-\eqref{eq:lasteq}, after replacing each matrix with its normalized (``tilde'') counterpart.

\section{Efficient Filtering Algorithms} \label{sec:filtering_algos}

The definitions of the filtered Loop and Star functions and of the quasi-Helmholtz Laplacian filters in Sections~\ref{sec:filtered_LS_bases} and \ref{sec:projector_filters} involve an SVD which, while ensuring a clear and compact theoretical treatment,  is in general computationally inefficient.
 This section will be devoted to presenting algorithms allowing for SVD-free matrix-vector products for the filtered  graph Laplacians  $\left(\mat\Sigma^T\mat\Sigma\right)_n$ and 
$\left(\mat\Lambda^T\mat\Lambda\right)_n$, which are the two key  operations on both approaches presented in the previous section. Moreover, while our treatment will deal with the graph matrices $\mat\Sigma$ and  $\mat\Lambda$, it is intended  that substantially the same strategies can be applied when replacing those matrices with their normalized counterparts $\tilde{\mat \Sigma}$ and $\tilde{\mat \Lambda}$ with minor modifications. In fact the additional products with the inverse square roots of (well-conditioned) Gram matrices can be obtained efficiently by using matrix function strategies \cite{higham_functions_2008}. 

\subsection{Power Method Filtering}

For filters with a filtering index that is independent on the total number of degrees of freedom, preconditioned inverse power methods \cite{golub_matrix_2012} yield the last singular vectors and singular values of the $\mat \Sigma^\T \mat \Sigma$ and $\mat \Lambda^\T \mat \Lambda$ matrices at the price of a constant number of matrix-vector products. Given that the matrices involved are sparse, the resulting method is linear in complexity and the filtered projectors can be efficiently obtained.
These schemes are well known and we do not provide extensive details here for the sake of brevity. We just mention that a special care should be put into using this schemes in the presence of degenerate spectra (arising from symmetries for example); under these conditions, the scheme presented in the following section should rather be preferred.

\subsection{Butterworth Matrix Filters}

An alternative strategy for the previous scenario, i.e. when a filter is needed that has filtering index which is independent on the total number of degrees of freedom, is provided by a matrix function and filtering approach. Given a scalar (squared) Butterworth filter of positive order $m$ and cutoff parameter $x_c > 0$, characterized by
\begin{equation}
    f_{m,x_c}(x) = \left(1 + (x/x_c)^{m}\right)^{-1}\,, \quad x \ge 0\,,
\end{equation}
the spectrum of a symmetric positive matrix $\mat A \in \mathbb{R}^{N\times N}$ composed of the set of singular values $\{\sigma_i(\mat A)\}_i$ can be filtered by generalizing $f_{m,x_c}$ to matrix arguments and applying it to $\mat A$, yielding the filtered matrix
\begin{equation}
    \mat A_\text{filt} \coloneqq f_{m,x_c}(\mat A) = \left(\matI + (\mat A/x_c)^{m}\right)^{-1}\,,
\end{equation}
with singular values $\{f_{m,x_c}\left(\sigma_i(\mat A)\right)\}_i$.
The filtered matrix $\left(\mat \Sigma^\T \mat \Sigma\right)_n$ can now be expressed as
\begin{equation}\label{eq:limit_proj}
    \left(\mat \Sigma^\T \mat \Sigma\right)_n = (\mat \Sigma^\T \mat \Sigma) \lim_{m\to\infty} f_{m,\sigma_n(\mat \Sigma^\T \mat \Sigma)} \left(\mat \Sigma^\T \mat \Sigma\right)\,.
\end{equation}
The presence of high exponents in \eqref{eq:limit_proj} may render its computation unstable. Hence we propose to use the following factorization formula that leverages the roots of unity
\begin{multline} \label{eq:factorized_butterworth}
    \left(\mat \Sigma^\T \mat \Sigma\right)_n = \left(\mat \Sigma^\T \mat \Sigma\right) \\\lim_{m\to\infty} \prod_{k=1}^m \left(\frac{\mat\Sigma^\T\mat\Sigma}{\sigma_n(\mat \Sigma^\T \mat \Sigma)} - \e^{\left(2 k + 1\right) \im\pi / N} \matI \right)^{-1}\,.
\end{multline}
For practical purposes the infinite products in this expression can be truncated at the desired precision. Regarding the value of
$\sigma_n(\mat \Sigma^\T \mat \Sigma)$, an approximation can be obtained either with ad-hoc heuristics or by the approximation $\sigma_n(\mat \Sigma^\T \mat \Sigma)\approx (N_s - n) / \| \left(\mat \Sigma^\T \mat \Sigma\right)^+ \|$. Finally, when the filtering point is a constant with respect to the number of unknowns, a multigrid approach is effective in providing the inverse required by \eqref{eq:factorized_butterworth}.
%\td{we need to double check the above statements}

\subsection{Filter Approximation via Chebyshev Polynomials}

When the filtering index is proportional to the number of unknowns, 
the computational burden of the two methods above can become high. In this regime we can leverage the ideas of polynomial preconditioning and graph wavelets \cite{johnsonPolynomialPreconditionersConjugate1983,ashbyComparisonAdaptiveChebyshev1992,hammond_wavelets_2011,levie2018cayleynets} and adopt a method based on a polynomial expansion of the spectral filter.
%\todo[inline,color=cyan]{FPA: This comments applies also to the introduction. We should find here a polynomial preconditioning approach that predates the wavelet guys already using such an approach.}

%\td{AM: Should we try to approximate heavyside directly? Leave this for a follow up paper}

Because we are interested in cases in which the filtering index is proportional to the number of degrees of freedom (for instance, $n = N_S / 2$) we can leverage a polynomial approximation of $f_{m,x_c}$ on the interval $[0, \sigma_{N_S}(\mat \Sigma^\T \mat \Sigma)]$; a natural basis for this approximation is that of the Chebyshev polynomials $\{T_n(x)\}_n$, defined by the recurrence relation
\begin{equation}
    T_n(x) = \begin{dcases}
        1 & \text{if } n = 0\\
        x & \text{if } n = 1\\
        2 x T_{n-1} (x) - T_{n-2} (x) & \text{otherwise.}
    \end{dcases}
\end{equation}
The approximated filtered matrix now reads
\begin{equation}\label{eq:cheby_filters_def}
    \left(\mat \Sigma^\T \mat \Sigma\right)_n \approx -\frac{c_0}{2} \matI + \sum_{k = 1}^{n_c} c_k T_k\left(\frac{\mat\Sigma^\T\mat\Sigma}{\sigma_n(\mat \Sigma^\T \mat \Sigma)}\right)\,,
\end{equation}
where the $c_n$ are the expansion coefficients of $f_{m,\sigma_n(\mat \Sigma^\T \mat \Sigma)}$ in the basis of the first $n_c+1$ Chebyshev polynomials. Algorithms for their computation can be found, among others, in \cite{press_numerical_2007}. Because the cutoff frequency of this filter is proportional to the number of unknowns and so is the domain size,
%\td{AM: I think that with another operator it would not work (eg any operator that does not have a linear spectrum in the number of unknowns.)}
the order of the polynomial that is required to obtain a given approximation of the Butterworth filter, will not need to be changed with increasing discretizations.
In other words, the filters obtained by following this approach will require the same number of sparse matrix-vector multiplication for increasing discretization when the filtering index will be proportional to the number of degrees of freedom.
It should be noted that in the transition region between the filters described in the previous two sections
(constant filtering index) and the scenario described here (filtering index will be proportional to the number of degrees of freedom) the Chebyshev approach decreases in efficiency and further treatments may be required \cite{levie2018cayleynets}. 

%\todo[color=yellow,inline]{We should add a comment about the complexity of each algorithm in each regime. Or at least add a comment and point to the Graph guys' papers}

\section{A First Application Case Scenario: Laplacian Filter Based Preconditioning} \label{sec:EFIE_application}
%\todo[inline]{CH: SHouldn't we use the normalized L and S matrices in this section? What is $U_{red}$} 

As a first application case scenario of the new filters introduced here, we will develop two families of preconditioners for the EFIE in \eqref{eq:disc_EFIE}. This equation is known to suffer from ill-conditioning both for decreasing frequency and average mesh length $h$ (phenomena known as the low-frequency and $h$-refinement breakdowns, respectively, see \cite{adrianElectromagneticIntegralEquations2021} and references therein). In the following we will cure both breakdowns by developing preconditioners based both on filtered functions decompositions and on quasi-Helmholtz Laplacian filters.

The reader should note that in this Section and in the subsequent ones, we will study the singular value spectrum of potentially singular matrices. When dealing with such matrices, the condition number will be defined as $\mathrm{cond}(\mat A) = \| \mat A \| \| \mat A^{+} \|$. Moreover, inverse powers of singular matrices in the following will always denote the corresponding positive power of the pseudoinverse of the matrix. 

\subsection{Filtered Bases Approach} \label{subsec:precond_filtered_LS}

The primal and dual Laplacians can be used to precondition the single layer and the hypersingular operator \cite{nedelec_acoustic_2001,mitharwal_multiplicative_2014,adrian_hierarchical_2017,oneilSecondkindIntegralEquations2018}, thus $\mat V_{\mat \Lambda}$, and $\mat V_{\mat \Sigma}$ followed by a diagonal preconditioning are valid bases for regularizing the vector and scalar potential parts of the EFIE.
In particular, for $\mat T_\mathrm{h}$, this results from the fact that an operator spectrally equivalent to the single layer can be obtained from $\mat T_\mathrm{h}$. In fact, noticing that $\mat T_\mathrm{h} = \mat \Sigma \mat R \mat \Sigma^\T$ \cite{zhao_integral_2000}, where $\mat R$ is the patch-function discretization of the single layer operator, i.e. $\left[ \mat R\right]_{mn} = \left<p_m, \op S p_n\right>$ with
\begin{equation}
    \left(\op S p\right) (\vt r) \coloneqq \int_\Gamma \frac{\e^{i k\| \vt{r} - \vt{r'} \|}}{4\pi\|\vt{r} - \vt{r'} \|} p(\vt{r'})  \dd S(\vt{r'})\,,
\end{equation}
and defining $\tmat R \coloneqq \mat G_{p}^{-1/2} \mat R \mat G_{p}^{-1/2}$, we obtain $\tmat T_\mathrm{h} = \tmat{\Sigma} \tmat R \tmat \Sigma^\T$. The equivalence between $\left(\tmat\Sigma^\T\tmat\Sigma\right)^+ \tmat \Sigma^\T \tmat T_\mathrm{h} \tmat \Sigma \left(\tmat\Sigma^\T\tmat\Sigma\right)^+$ and $\tmat R$ thus follows.
To conclude the reasoning, we note that, because
\begin{multline}
    \left(\tmat\Sigma^\T\tmat\Sigma\right)^{1/4}\tmat R\left(\tmat\Sigma^\T\tmat\Sigma\right)^{1/4}\\= \tmat V_{\tmat \Sigma} \left(\tmat S_{\tmat\Sigma}^\T\tmat S_{\tmat\Sigma}\right)^{1/4} \tmat V_{\tmat \Sigma}^\T \tmat R \tmat V_{\tmat \Sigma} \left(\tmat S_{\tmat\Sigma}^\T\tmat S_{\tmat\Sigma}\right)^{1/4} \tmat V_{\tmat \Sigma}^\T\,,
\end{multline}
is well conditioned for increasing discretization---as a consequence of the results proven in \cite{oneilSecondkindIntegralEquations2018}, since $\tmat\Sigma^\T\tmat\Sigma$ is a valid discretization of a Laplacian matrix \cite{arnold2007compatible}---we have
\begin{multline} \label{eq:standard_laplacian_precond_Th}
    \mathrm{cond}\left(\tmat V_{\tmat \Sigma} \left(\tmat S_{\tmat\Sigma}^\T\tmat S_{\tmat\Sigma}\right)^{1/4} \tmat V_{\tmat \Sigma}^\T \left(\tmat\Sigma^\T\tmat\Sigma\right)^+ \tmat \Sigma^\T \tmat T_\mathrm{h}\right.\\\left.\tmat \Sigma \left(\tmat\Sigma^\T\tmat\Sigma\right)^+ \tmat V_{\tmat \Sigma} \left(\tmat S_{\tmat\Sigma}^\T \tmat S_{\tmat \Sigma}\right)^{1/4} \tmat V_{\tmat \Sigma}^\T \right)=O(1)\,, h\to0\,.
\end{multline}
%\td{Each diagonal matrix should have the strength of $\left(\frac{1}{h}\right)^{\frac{1}{2}}$}
The reader should note that, since $\tmat V_{\tmat \Sigma}$ is unitary, we also have
\begin{multline}
    \mathrm{cond}\left(\tmat V_{\tmat \Sigma} \left(\tmat S_{\tmat\Sigma}^\T\tmat S_{\tmat\Sigma}\right)^{1/4} \tmat V_{\tmat \Sigma}^\T \left(\tmat\Sigma^\T\tmat\Sigma\right)^+ \tmat \Sigma^\T \tmat T_\mathrm{h}\right.\\\left.\tmat \Sigma \left(\tmat\Sigma^\T\tmat\Sigma\right)^+ \tmat V_{\tmat \Sigma} \left(\tmat S_{\tmat\Sigma}^\T \tmat S_{\tmat \Sigma}\right)^{1/4} \tmat V_{\tmat \Sigma}^\T \right) =\\
    \mathrm{cond}\left(\left(\tmat S_{\tmat\Sigma}^\T\tmat S_{\tmat\Sigma}\right)^{1/4} \tmat V_{\tmat \Sigma}^\T \left(\tmat\Sigma^\T\tmat\Sigma\right)^+ \tmat \Sigma^\T \tmat T_\mathrm{h}\right.\\\left.\tmat \Sigma \left(\tmat\Sigma^\T\tmat\Sigma\right)^+ \tmat V_{\tmat \Sigma} \left(\tmat S_{\tmat\Sigma}^\T \tmat S_{\tmat \Sigma}\right)^{1/4}\right)\,.
\end{multline}
Such an approach would require the computation of the matrix $\tmat V_{\tmat \Sigma}$ and $\tmat S_{\tmat \Sigma}$ which are prohibitively expensive to obtain. A key observation, however, is that we do not need to use the entire diagonal of $\tmat S_{\tmat\Sigma}$, but a logarithmic sampling of it will suffice.
In other words, define $\vec D_{\tmat\Sigma}$ the vector containing the entries of the diagonal of $\tmat S_{\tmat\Sigma}^\T\tmat S_{\tmat\Sigma}$ and define the block diagonal matrix
\begin{multline}
    \tmat D_{\tmat\Sigma,\alpha}=\mathrm{diag}\left(\left[\vec D_{\tmat\Sigma}\right]_{N_S - N_{S,\alpha} + 1} \mat I_{N^\mathrm{rem}_{S,\alpha}},\right.\\\left.\left[\vec D_{\tmat\Sigma}\right]_{N_S - \frac{N_{S,\alpha}}{\alpha} + 1} \mat I_{\frac{N_{S,\alpha}}{\alpha}},
    \ldots,\left[\vec D_{\tmat\Sigma}\right]_{N_S}\mat I_1 \right)\,,
\end{multline}
where $N_{S,\alpha} = \alpha^{\lfloor \log_\alpha(N_S) \rfloor}$, $N^\mathrm{rem}_{S,\alpha} = N_S - \left(1 - N_{S,\alpha}\right)\left(1-\alpha\right)^{-1}$, and $\mat I_n$ is the identity matrix of size $n$, or, more programmatically,
\begin{equation}
    \left[\tmat D_{\tmat\Sigma,\alpha}\right]_{ii} = \left[\vec D_{\tmat\Sigma}\right]_{f_{\tmat \Sigma}(i)}\,,
\end{equation}
with  $f_{\tmat \Sigma}(i)=  N_S - \alpha^{\lfloor \log_\alpha(N_S - i + 1) \rfloor} + 1$.
%the complexity in indexing is introduced here to remain consistent with the standard SVD definition, but can be avoided altogether by defining an SVD in which the singular values are in increasing order.
Note that the construction of this matrix only requires explicit knowledge of $\log_\alpha(N_S)$ terms of $\vec D_{\tmat \Sigma}$.
%This matrix, which contains $\log_\alpha(N_{S,\alpha}) + 1$ blocks, is depicted in Fig.~\ref{fig:diag}, for clarity. Few passages---omitted here---suffice to show that
Few passages---omitted here---suffice to show that
\begin{multline}
    {\rm cond}\left(\tmat D^{1/4}_{\tmat\Sigma,\alpha}\tmat V_{\tmat \Sigma}^\T\left(\tmat\Sigma^\T\tmat\Sigma\right)^+\tmat\Sigma^\T\tmat T_\mathrm{h}\tmat\Sigma\left(\tmat\Sigma^\T\tmat\Sigma\right)^+\tmat V_{\tmat \Sigma}\tmat D^{1/4}_{\tmat\Sigma,\alpha}\right)\\=O(\alpha)=O(1)\,, h\to0\,,
\end{multline}
which is reminiscent of hierarchical strategies (see \cite{adrianElectromagneticIntegralEquations2021} and references therein). Because $\tmat V_{\tmat \Sigma}$ is unitary, we obtain equivalently
\begin{multline}\label{eqnolabel1}
  {\rm cond}\left(\tmat V_{\tmat \Sigma}\tmat D^{1/4}_{\tmat\Sigma,\alpha} \tmat V_{\tmat \Sigma}^\T\left(\tmat\Sigma^\T\tmat\Sigma\right)^+\tmat\Sigma^\T\tmat T_\mathrm{h}\right.\\\left.\tmat\Sigma\left(\tmat\Sigma^\T\tmat\Sigma\right)^+\tmat V_{\tmat \Sigma}\tmat D^{1/4}_{\tmat\Sigma,\alpha}\tmat V_{\tmat \Sigma}^\T\right)=O(\alpha)=O(1)\,.
\end{multline}
%\td{AM: I chose to introduce the additional laplacian only here so that we do not look retarded as it will simplify everywhere else...}
%\todo[color=green,inline]{FPA: just state you want the basis to appear, for FPA to do.}
%While this preconditioning strategy is sound in itself, it can be further modified to leverage the filtered basis introduce in Section~\ref{} by introducing an additional graph Laplacian in \eqref{eq:standard_laplacian_precond_Th} and adjusting the exponent of $\tmat S_{\tmat \Sigma}^T \tmat S_{\tmat \Sigma}$. In particular, we have,

This preconditioning strategy can be slightly altered to leverage the filtered basis presented in Section~\ref{sec:filtered_LS_bases} by introducing an additional Laplacian in \eqref{eq:standard_laplacian_precond_Th} and adjusting the exponent of $\tmat S_{\tmat \Sigma}^\T \tmat S_{\tmat \Sigma}$ accordingly. In particular, we have
\begin{multline} \label{eq:filt_star_add_lalpacian}
    \tmat\Sigma \left(\tmat\Sigma^\T\tmat\Sigma\right)^+ \tmat V_{\tmat \Sigma} \left(\tmat S_{\tmat \Sigma}^\T \tmat S_{\tmat \Sigma} \right)^{1/4} \tmat V_{\tmat \Sigma}^\T  =\\ \tmat\Sigma \left(\tmat\Sigma^\T\tmat\Sigma\right)^+ \left(\tmat\Sigma^\T\tmat\Sigma\right) \tmat V_{\tmat \Sigma} \left(\tmat S_{\tmat \Sigma}^\T \tmat S_{\tmat \Sigma} \right)^{-3/4} \tmat V_{\tmat \Sigma}^\T\,,
\end{multline}
which, following the reasoning detailed above, means that
\begin{equation}
\tmat B_{\tmat \Sigma} \coloneqq \tmat\Sigma \left(\tmat\Sigma^\T\tmat\Sigma\right)^+ \left(\tmat\Sigma^\T\tmat\Sigma\right) \tmat V_{\tmat \Sigma}\tmat D^{-3/4}_{\tmat\Sigma,\alpha}\tmat V_{\tmat \Sigma}^\T
\end{equation}
is a valid left and right symmetric preconditioner for $\tmat T_\mathrm{h}$. Finally, thanks to the properties introduced in Section~\ref{sec:filtered_LS_bases}, we have
\begin{multline} \label{eqnolabel2}
\tmat\Sigma\left(\tmat\Sigma^\T\tmat\Sigma\right)^+\left(\tmat\Sigma^\T\tmat\Sigma\right)\tmat V_{\tmat \Sigma}\tmat D^{-3/4}_{\tmat\Sigma,\alpha}\tmat V_{\tmat \Sigma}^\T=\\ \sum_{l=2}^{N_{S,\alpha}} \left({\tmat\Sigma}_{\alpha^{l}-1}-{\tmat\Sigma}_{\alpha^{l-1}-1}\right) \left[\vec D_{\tmat \Sigma}\right]_{N_S - \alpha^{l-1}+1}^{-3/4}\\
+ \left({\tmat\Sigma}-{\tmat\Sigma}_{\alpha^{N_{S,\alpha}}-1}\right) \left[\vec D_{\tmat \Sigma}\right]_{N_S - N_{S,\alpha}+1}^{-3/4}
\eqqcolon{\tmat\Sigma}_{\mathrm{p},\alpha}
\end{multline} 
and thus from  \eqref{eq:filt_star_add_lalpacian} and \eqref{eqnolabel2} it follows that
\begin{equation} \label{eq:precond_LS_Th}
       {\rm cond}\left( \tmat\Sigma_{\mathrm{p},\alpha}^\T \tmat T_\mathrm{h}  \tmat\Sigma_{\mathrm{p},\alpha}\right)=O(1)\,, h\to 0\,.
\end{equation}
%\td{Here instead the alpha the strength of $\left(\frac{1}{h}\right)^{-\frac{3}{2}}$ because of the extra $\left(\tmat\Sigma^\T\tmat\Sigma_n\right)$ in the definition}

A similar reasoning for $\tmat T_\mathrm{s}$, following from the preconditioning of the hypersingular operator, leads to
\begin{equation}
    {\rm cond}\left(\left(\tmat S_{\tmat\Lambda}^\T\tmat S_{\tmat\Lambda}\right)^{-1/4}\tmat V_{\tmat \Lambda}^\T\tmat\Lambda^\T\tmat T_\mathrm{s} \tmat\Lambda\tmat V_{\tmat \Lambda}\left(\tmat S_{\tmat\Lambda}^\T\tmat S_{\tmat\Lambda}\right)^{-1/4}\right)=O(1)
\end{equation}
and
\begin{multline}
    {\rm cond}\left(\tmat V_{\tmat \Lambda}\left(\tmat S_{\tmat\Lambda}^\T\tmat S_{\tmat\Lambda}\right)^{-1/4}\tmat V_{\tmat \Lambda}^\T\tmat\Lambda^\T\tmat T_\mathrm{s}\right.\\\left. \tmat\Lambda\tmat V_{\tmat \Lambda}\left(\tmat S_{\tmat\Lambda}^\T\tmat S_{\tmat\Lambda}\right)^{-1/4}\tmat V_{\tmat \Lambda}^\T\right)=O(1)\,, h \to 0\,.
\end{multline}
%\td{Each diagonal matrix should have the strength of $\left(\frac{1}{h}\right)^{-\frac{1}{2}}$}
From this, following a dual reasoning as the one of the previous section, we obtain
\begin{equation}
       {\rm cond}\left( {\tmat\Lambda}_{\mathrm{p},\alpha}^\T \tmat T_\mathrm{s}   {\tmat\Lambda}_{\mathrm{p},\alpha}\right)=O(1)\,, h \to 0\,.
\end{equation}
where
\begin{multline}
\tmat\Lambda\left(\tmat\Lambda^\T\tmat\Lambda\right)^+\left(\tmat\Lambda^\T\tmat\Lambda\right)\tmat V_{\tmat \Lambda}\tmat D^{-1/4}_{\tmat\Lambda,\alpha}\tmat V_{\tmat \Lambda}^\T=\\ \sum_{l=2}^{N_{L,\alpha}} \left({\tmat\Lambda}_{\alpha^{l}-1}-{\tmat\Lambda}_{\alpha^{l-1}-1}\right) \left[\vec D_{\tmat \Lambda}\right]_{N_L - \alpha^{l-1}+1}^{-1/4}\\
+ \left({\tmat\Lambda}-{\tmat\Lambda}_{\alpha^{N_{L,\alpha}}-1}\right) \left[\vec D_{\tmat \Lambda}\right]_{N_L - N_{L,\alpha}+1}^{-1/4}
\eqqcolon \tmat\Lambda_{\mathrm{p},\alpha}\,,
\end{multline}
and
\begin{multline}
    \tmat D_{\tmat\Lambda,\alpha}=\mathrm{diag}\left(\left[\vec D_{\tmat\Lambda}\right]_{N_L - N_{L,\alpha} + 1} \mat I_{N^\mathrm{rem}_{L,\alpha}},\right.\\\left.\left[\vec D_{\tmat\Lambda}\right]_{N_L - \frac{N_{L,\alpha}}{\alpha} + 1} \mat I_{\frac{N_{L,\alpha}}{\alpha}},
    \ldots,\left[\vec D_{\tmat\Lambda}\right]_{N_L}\mat I_1 \right)\,,
\end{multline}
with $N^\mathrm{rem}_{L,\alpha} = N_L - \left(1 - N_{L,\alpha}\right)\left(1-\alpha\right)^{-1}$, $\vec D_{\tmat\Lambda}$ the vector containing the elements of the diagonal of $\tmat S_{\tmat \Lambda}^\T\tmat S_{\tmat \Lambda}$, and $N_{L,\alpha} = \alpha^{\lfloor \log_\alpha(N_L) \rfloor}$.
%\td{Here the coefficients alpha are confirmed to have the strength of $\left(\frac{1}{h}\right)^{-\frac{1}{2}}$}

%\td{AM: Check if the norm scaling is still necessary after the last changes.}

The previous preconditioners can then be combined to obtain a complete regularization of the EFIE system, for both low-frequency and h-refinement breakdowns, that reads
\begin{equation}\label{eq:EFIE_wavelet}
\tmat W^\T \tmat T \tmat W \vec{\tilde{j}} = \tmat W^\T \vec{\tilde{v}}\,,
\end{equation}
where $\vec{\tilde{v}} = \mat G^{-1/2} \vec v$, $\vec j= \mat G^{-1/2} \tmat W \vec{\tilde{j}}$, $\tmat W = \begin{bmatrix} \sqrt{c_{\tmat\Lambda}}{\tmat\Lambda}_{\mathrm{p},\alpha} & \sqrt{c_{\tmat\Sigma}}{\tmat\Sigma}_{\mathrm{p},\alpha}
\end{bmatrix}$, $c_{\tmat\Sigma} = \|\tmat\Sigma^\T_{\mathrm{p},\alpha} \tmat T_\mathrm{h} \tmat\Sigma_{\mathrm{p},\alpha}\|^{-1}$, $c_{\tmat\Lambda} = \|\tmat\Lambda^{\T}_{\mathrm{p},\alpha} \tmat T_\mathrm{s}  \tmat\Lambda_{\mathrm{p},\alpha}\|^{-1}$, and where we assume that the appropriate number of columns have been removed from $\tmat\Sigma_{\mathrm{p},\alpha}$ and $\tmat\Lambda_{\mathrm{p},\alpha}$ (e.g. 1 column must be removed from each for a simply connected, closed scatterer) to account for the linear dependence in the underlying Loop and Star bases \cite{wilton_topological_1983}, as is done in standard Loop-Star preconditioning. The reader should note that, as in the case of standard Loop-Star functions, this operations will create a small number of isolated singular values, that however will not impact the convergence properties of the preconditioned equation. This effect will not be present instead in the scheme of next Section. 
The $h$-refinement regularization effect of this preconditioner can be deduced from the previous derivations for each of the potentials \cite{boubendir_well-conditioned_2014}. 
\begin{comment}
The h-refinement regularization effect of this preconditioner can be deduced from the previous derivations since the left-hand-side of \eqref{eq:EFIE_wavelet} can be written the sum of a well conditioned matrix and one arising from the discretization of compact operators
\begin{align}
\begin{split}
    &\mat W^T \mat{T} \mat W\\
    {}&=\begin{bmatrix}
        \frac{1}{c_{\mat\Lambda}}{\mat\Lambda}^{T,\text{trunc}}_{h,\alpha} \mat{Ts} {\mat\Lambda}^{\text{trunc}}_{h,\alpha} & \matO \\
        \matO & \frac{1}{c_{\mat\Sigma}} {\mat\Sigma}^{T,\text{trunc}}_{h,\alpha} \mat T_\mathrm{h} {\mat\Sigma}^{\text{trunc}}_{h,\alpha}
    \end{bmatrix}\\
    &+\begin{bmatrix} \matO & \frac{1}{\sqrt{c_{\mat\Lambda}c_{\mat\Sigma}}}{\mat\Lambda}^{T,\text{trunc}}_{h,\alpha} \mat T_\mathrm{s}  {\mat\Sigma}^{\text{trunc}}_{h,\alpha} \\
    \frac{1}{\sqrt{c_{\mat\Lambda}c_{\mat\Sigma}}}{\mat\Sigma}^{T,\text{trunc}}_{h,\alpha} \mat T_\mathrm{s}  {\mat\Lambda}^{\text{trunc}}_{h,\alpha}& \frac{1}{c_{\mat\Sigma}}{\mat\Sigma}^{T,\text{trunc}}_{h,\alpha} \mat T_\mathrm{s}  {\mat\Sigma}^{\text{trunc}}_{h,\alpha}
    \end{bmatrix}\,.
    \end{split}
\end{align}
\end{comment}
The low frequency regularization, can be demonstrated following the same reasoning as for standard Loop-Star approaches \cite{adrianElectromagneticIntegralEquations2021}, since the new filtered bases retain the crucial properties that made Loop-Star so adapted low-frequency regularization in the first place---${\tmat\Lambda}^{\T}_{\mathrm{p},\alpha} \tmat T_\mathrm{h} = \matO$, $\tmat T_\mathrm{h} {\tmat\Lambda}_{\mathrm{p},\alpha} = \matO$, and ${\tmat\Lambda}^{\T}_{\mathrm{p},\alpha} {\tmat\Sigma}_{\mathrm{p},\alpha} = \matO$.
Finally, we have
\begin{equation}
   \mathrm{cond}\left(\tmat W^\T \tmat T \tmat W\right)=O(1)\,, \text{ when } h\to 0\,,k \to 0\,.
\end{equation}

\subsection{Quasi-Helmholtz Filters Approach} \label{subsec:precond_filtered_projectors}

In several application scenarios, an explicit quasi-Helmholtz decomposition, such as the Loop-Star decomposition, is not necessary, and quasi-Helmholtz projectors \cite{adrianElectromagneticIntegralEquations2021} could be used instead. Similarly, instead of using filtered Loop-Star preconditioning approaches, basis-free approaches, based on the quasi-Helmholtz Laplacian filters, will often be more effective. This Section will explore this approach that, as an additional advantage, will also have the avoidance of the burden  of global-loop detection for multiply connected scatterers.
%It is also, as will be shown in Section~\ref{sec:filtering_algos}, amenable to efficient computation techniques, that elude the loop-star approach.

Following the same philosophy as in Section~\ref{subsec:precond_filtered_LS}, we will form preconditioners for the solenoidal part of $\tmat T_\mathrm{s} $ and for $\tmat T_\mathrm{h}$ that will then be combined into a full EFIE preconditioner. We can transition from a basis-based Helmholtz decomposition to a projector based Helmholtz decomposition by leveraging the correspondences between $\tmat\Sigma$ and $\tmat\Lambda$ and their respective projectors $\tmat P^\Sigma$ and $\tmat P^\Lambda$.
%Similarly, a quasi-Helmholtz Filters based preconditioner can be obtained from after replacing bla with bla
%In the same way that we added an additional Laplacian in \eqref{eq:filt_star_add_lalpacian} to form a preconditioner based on the filtered basis, we can add an additional discretized divergence operator to the right of the preconditioner and adjust the exponent of $\tmat D_{\tmat \Sigma,\alpha}$ to obtain a preconditioner based on filters. 
In particular, because $\tmat B_{\tmat \Sigma}$
%\begin{equation}
 %\coloneqq \tmat\Sigma\left(\tmat\Sigma^\T\tmat\Sigma\right)^+\left(\tmat\Sigma^\T\tmat\Sigma\right)\tmat V_{\tmat \Sigma}\tmat D^{-3/4}_{\tmat\Sigma}
%\end{equation}
was a valid preconditioner for $\tmat T_\mathrm{h}$ (equations \eqref{eqnolabel2} and \eqref{eq:precond_LS_Th}), $\begin{bmatrix} \tmat B_{\tmat \Sigma} & \matO \end{bmatrix}$, once applied left and right to $\tmat T_\mathrm{h}$ will yield a block diagonal matrix which is well conditioned away from its large nullspace. This, in turns, means that $\begin{bmatrix} \tmat C_{\tmat \Sigma} & \matO \end{bmatrix}$, with $\tmat C_{\tmat \Sigma}=\tmat\Sigma\left(\tmat\Sigma^\T\tmat\Sigma\right)^+\left(\tmat\Sigma^\T\tmat\Sigma\right)\tmat V_{\tmat \Sigma}\tmat D^{-5/4}_{\tmat\Sigma} \tmat V_{\tmat \Sigma}^\T \tmat V_{\tmat \Sigma} \tmat D^{1/2}_{\tmat\Sigma}$, will also yield a well-conditioned (up to its nullspace) matrix. Finally, because multiplications by unitary matrices do not compromise conditioning properties, we can form the preconditioner
\begin{multline}
    \tmat\Sigma\left(\tmat\Sigma^\T\tmat\Sigma\right)^+\left(\tmat\Sigma^\T\tmat\Sigma\right) \tmat V_{\tmat \Sigma}\tmat D^{-5/4}_{\tmat\Sigma} \tmat V_{\tmat \Sigma}^\T \tmat V_{\tmat \Sigma} \begin{bmatrix}\tmat D^{1/2}_{\tmat\Sigma} & \matO \end{bmatrix} \tmat U_{\tmat \Sigma}^\T =\\
    \tmat\Sigma\left(\tmat\Sigma^\T\tmat\Sigma\right)^+\tmat V_{\tmat \Sigma}\tmat D^{-1/4}_{\tmat\Sigma} \tmat V_{\tmat \Sigma}^\T \tmat \Sigma^\T\,.
\end{multline}
This allows us to form the preconditioner $\tmat Q^{\tmat \Sigma}_{\mathrm{p},\alpha}$, of additive Schwarz type, based on quasi-Helmholtz filters
\begin{multline}
\tmat\Sigma\left(\tmat\Sigma^\T\tmat\Sigma\right)^+\tmat V_{\tmat \Sigma}\tmat D^{-1/4}_{\tmat\Sigma} \tmat V_{\tmat \Sigma}^\T \tmat \Sigma^\T=\\ \sum_{l=2}^{N_{S,\alpha}} \left(\tmat P^{\tmat \Sigma}_{\alpha^{l}-1}-\tmat P^{\tmat\Sigma}_{\alpha^{l-1}-1}\right) \left[\vec D_{\tmat \Sigma}\right]_{N_S - \alpha^{l-1}+1}^{-1/4}\\
+ \left(\tmat P^{\tmat \Sigma}-\tmat P^{\tmat\Sigma}_{\alpha^{N_{S,\alpha}}-1}\right) \left[\vec D_{\tmat \Sigma}\right]_{N_S - N_{S,\alpha}+1}^{-1/4} \eqqcolon \tmat Q^{\tmat \Sigma}_{\mathrm{p},\alpha}
\end{multline}
for which
\begin{equation}\label{eq:defSigmatilde}
       {\rm cond}\left( \tmat Q^{\tmat\Sigma}_{\mathrm{p},\alpha} \tmat T_\mathrm{h} \tmat Q^{\tmat\Sigma}_{\mathrm{p},\alpha}\right)=O(1)\,, h\to 0\,.
\end{equation}

Similarly, a preconditioner for the solenoidal part of $\tmat T_\mathrm{s} $ is
\begin{multline}
\tmat Q^{\tmat \Lambda}_{\mathrm{p},\alpha} \coloneqq \sum_{l=2}^{N_{L,\alpha}} \left(\tmat P^{\tmat \Lambda}_{\alpha^{l}-1}-\tmat P^{\tmat\Lambda}_{\alpha^{l-1}-1}\right) \left[\vec D_{\tmat \Lambda}\right]_{N_L - \alpha^{l-1}+1}^{1/4}\\
+ \left(\tmat P^{\tmat \Lambda}-\tmat P^{\tmat\Lambda}_{\alpha^{N_{L,\alpha}}-1}\right) \left[\vec D_{\tmat \Lambda}\right]_{N_L - N_{L,\alpha}+1}^{1/4}
\end{multline}
for which
\begin{equation}\label{eq:defLambdatilde}
       {\rm cond}\left(\tmat Q^{\tmat\Lambda}_{\mathrm{p},\alpha} \tmat T_\mathrm{s}  \tmat Q^{\tmat\Lambda}_{\mathrm{p},\alpha}\right)=O(1)\,, h\to0\,.
\end{equation}

The full EFIE preconditioner is then an appropriate linear combination of the solenoidal and non-solenoidal preconditioners above to cure also the low-frequency breakdown. In particular we define
\begin{equation}
    \tmat Q = \sqrt{b_{\tmat \Lambda}} \tmat Q^{\tmat \Lambda}_{\mathrm{p},\alpha} + i\sqrt{b_{\tmat \Sigma}} \tmat Q^{\tmat \Sigma}_{\mathrm{p},\alpha} + \sqrt{b_{\tmat H}} \tmat P^{H}\,,
\end{equation}
where $\tmat P^{H}=\mat I-\tilde{\mat P}^{\Sigma}-\tilde{\mat P}^\Lambda$
and
\begin{align}
    b_{\tmat \Lambda} &= \| \tmat Q^{\tmat \Lambda}_{\mathrm{p},\alpha} \tmat T_\mathrm{s}  \tmat Q^{\tmat \Lambda}_{\mathrm{p},\alpha}\|^{-1}\,,\\
    b_{\tmat \Sigma} &= \| \tmat Q^{\tmat \Sigma}_{\mathrm{p},\alpha} \tmat T_\mathrm{h} \tmat Q^{\tmat \Sigma}_{\mathrm{p},\alpha}\|^{-1}\,,\\
    b_{\tmat H} &= \| \tmat P^{\tmat H} \tmat T_\mathrm{s}  \tmat P^{H}\|^{-1}\,,
\end{align}
account for the frequency-scaling of the operators and the diameter of $\Gamma$.
The preconditioned EFIE system is
\begin{equation}\label{eq:EFIE_proj}
    \tmat Q \tmat{T} \tmat Q \tilde{\vec j}_\mathrm{qH} = \tmat Q \vec{\tilde{v}}\,,
\end{equation}
with $\vec j = \mat G^{-1/2} \tmat Q \tilde{\vec j}_\mathrm{qH}$.

\section{Implementation related details and further improvements} \label{eq:implementation_details}

In addition to the efficient filtering algorithms presented in Section~\ref{sec:filtering_algos}, obtaining a fast and efficient implementation of the  proposed preconditioning scheme based on filtered projectors requires that particular attention be given to parts of their implementation. First, all the terms of the form $\mat T_\mathrm{h} \mat Q^{\mat \Lambda}_{\mathrm{p},\alpha}$, $\mat Q^{\mat \Lambda}_{\mathrm{p},\alpha} \mat T_\mathrm{h}$, $\mat P^H \mat T_\mathrm{h}$, or $\mat T_\mathrm{h} \mat P^H$ must be explicitly set to $\matO$ to avoid numerical instabilities. Further treatments on the right hand side and on the solution vector, are required to ensure that the solution of the system remains accurate until arbitrarily low frequencies. These treatments are straightforward generalization of those required for standard quasi-Helmholtz preconditioning techniques that can be found in \cite{adrianElectromagneticIntegralEquations2021}.

%To ensure an effective dense discretization preconditioning effect, the schemes introduced in Section~\ref{subsec:precond_filtered_LS} and Section~\ref{subsec:precond_filtered_projectors} must be slightly modified. Upon discretization, in the case of general and potentially non-uniformly discretized geometries, the diagonal preconditioning based on the theoretical Laplacian eigenvalues provides a solid theoretical approach but sub-optimal conditioning results. A more resilient approach instead is to use operator norms, the new preconditioners then become
The condition numbers obtained when employing the schemes introduced in Section~\ref{subsec:precond_filtered_LS} and Section~\ref{subsec:precond_filtered_projectors}, while stable, can be further brought down by slightly modifying the preconditioners. The diagonal preconditioning based on the theoretical Laplacian eigenvalues can be altered to instead employ matrix norms; the new preconditioners then become
\begin{multline} \label{eq:norm_star_proj_scaling}
\mat Q^{\mat \Sigma}_{\mathrm{p},\alpha} = \sum_{l=2}^{N_{S,\alpha}} \left(\mat P^{\mat \Sigma}_{\alpha^{l}-1}-\mat P^{\mat\Sigma}_{\alpha^{l-1}-1}\right) b_l
+\\\left(\mat P^{\mat \Sigma}-\mat P^{\mat\Sigma}_{\alpha^{N_{S,\alpha}}-1}\right) b_{N_{S,\alpha}+1}\,,
\end{multline}
where
\begin{gather}
\begin{split}
    b_l = \left\| \left(\mat P^{\mat \Sigma}_{\alpha^{l}-1}-\mat P^{\mat\Sigma}_{\alpha^{l-1}-1}\right)^\T \mat T_\mathrm{h} \left(\mat P^{\mat \Sigma}_{\alpha^{l}-1}-\mat P^{\mat\Sigma}_{\alpha^{l-1}-1}\right)\right\|^{-1/2}\,,\\2 \le l \le N_{S,\alpha}\,,
    \end{split}\\
    b_{N_{S,\alpha}+1} = \left\| \left(\mat P^{\mat \Sigma}-\mat P^{\mat\Sigma}_{\alpha^{N_{S,\alpha}}-1}\right)^\T \mat T_\mathrm{h} \left(\mat P^{\mat \Sigma}-\mat P^{\mat\Sigma}_{\alpha^{N_{S,\alpha}}-1}\right)\right\|^{-1/2}\,.
\end{gather}
The same modification can be performed for $\mat Q^{\mat\Lambda}_{\mathrm{p},\alpha}$ that becomes
\begin{multline} \label{eq:norm_loop_proj_scaling}
\mat Q^{\mat \Lambda}_{\mathrm{p},\alpha} = \sum_{l=2}^{N_{L,\alpha}} \left(\mat P^{\mat \Lambda}_{\alpha^{l}-1}-\mat P^{\mat\Lambda}_{\alpha^{l-1}-1}\right) d_l
+\\\left(\mat P^{\mat \Lambda}-\mat P^{\mat\Lambda}_{\alpha^{N_{L,\alpha}}-1}\right) d_{N_{S,\alpha}+1}\,,
\end{multline}
with
\begin{gather}
\begin{split}
    d_l = \left\| \left(\mat P^{\mat \Lambda}_{\alpha^{l}-1}-\mat P^{\mat\Lambda}_{\alpha^{l-1}-1}\right)^\T \mat T_\mathrm{s}  \left(\mat P^{\mat \Lambda}_{\alpha^{l}-1}-\mat P^{\mat\Lambda}_{\alpha^{l-1}-1}\right)\right\|^{-1/2}\,,\\2 \le l \le N_{L,\alpha}\,,
    \end{split}\\
    d_{N_{S,\alpha}+1} = \left\| \left(\mat P^{\mat \Lambda}-\mat P^{\mat\Lambda}_{\alpha^{N_{L,\alpha}}-1}\right)^\T \mat T_\mathrm{s}  \left(\mat P^{\mat \Lambda}-\mat P^{\mat\Lambda}_{\alpha^{N_{L,\alpha}}-1}\right)\right\|^{-1/2}\,.
\end{gather}
To ensure that the overall complexity of the algorithm is not increased, the values of $\{b_l\}_l$ and $\{d_l\}_l$ can be efficiently computed using, for example, power methods.
The reader should note that the preconditioning approach delineated above requires filter profiles with support both proportional to and independent from the number of unknowns, which can be efficiently obtained with the approaches described in Section~\ref{sec:filtering_algos}. As said in the previous Section, filters in the transition region could be less efficient to obtain, as the Chebyshev approach decreases in efficiency away from the middle of the spectrum \cite{levie2018cayleynets}. All preconditioning real case scenarios presented here, however, are not  impacted by this fact as shown in Section~\ref{sec:numerical_results}. 

\section{Numerical Results} \label{sec:numerical_results}

\begin{figure}
    \centering
    \includegraphics{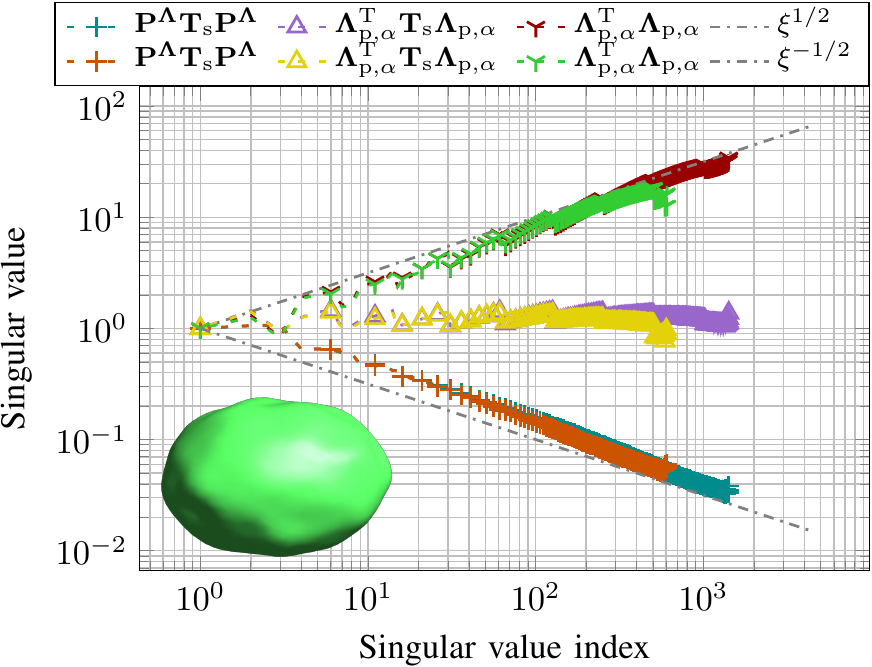}
    \caption{Spectrum of the solenoidal part of the vector potential, its preconditioner, and its preconditioned counterpart. These spectra have been obtained for a smoothly-deformed sphere with a maximum diameter of \SI{7.17}{\meter} (see insert), a frequency of \SI{e6}{\hertz}, and for two different average edge lengths \SI{0.31}{\meter} and \SI{0.20}{\meter}. The spectra have been normalized so that their first singular value is one, for readability. Perfect filters built out of SVD have been used in these results.}
    \label{fig:spectrum_Ts_filtered_loops}
\end{figure}

\begin{figure}
    \centering
    \includegraphics{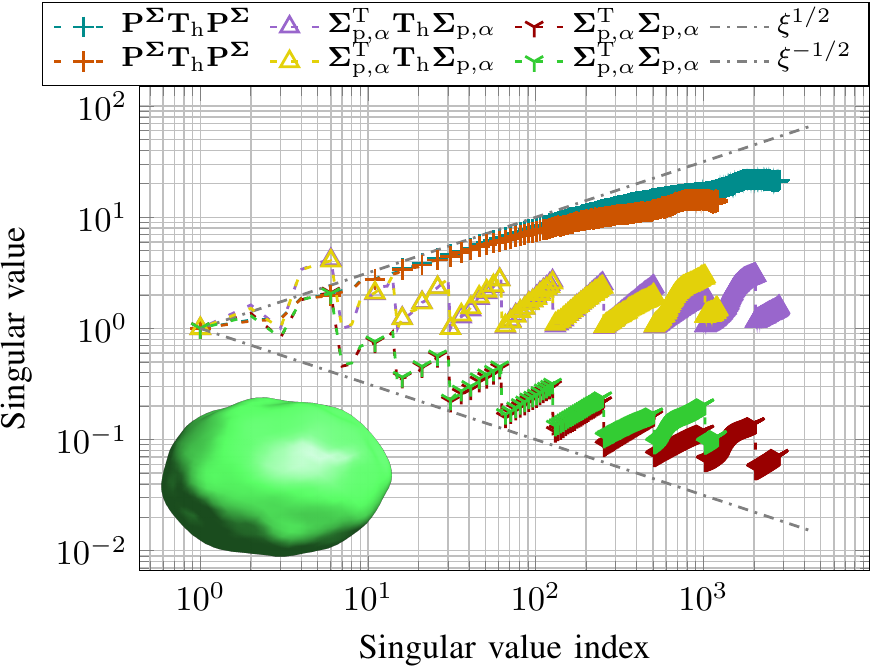}
    \caption{Spectrum of the non-solenoidal part of the scalar potential, its preconditioner, and its preconditioned counterpart. These spectra have been obtained for a smoothly-deformed sphere with a maximum diameter of \SI{7.17}{\meter} (see insert), a frequency of \SI{e6}{\hertz}, and for two different average edge lengths \SI{0.31}{\meter} and \SI{0.20}{\meter}. The spectra have been normalized so that their first singular value is one, for readability. Perfect filters built out of SVD have been used in these results.}
    \label{fig:spectrum_Th_filtered_stars}
\end{figure}

All  numerical results presented in this section have been obtained with non-normalized matrices ($\mat \Lambda$, $\mat \Sigma$) to illustrate that graph matrices are often enough for practical cases. Equally good or superior performance, however, can be obtained by using   normalized  matrices ($\tilde{\mat \Lambda}$, $\tilde{\mat \Sigma}$) instead. In the first set of examples we have leveraged perfect filters obtained by SVD before presenting results based on SVD-free approaches.
The filtered Loop-Star preconditioning approach presented in Section~\ref{subsec:precond_filtered_LS} leverages the spectral equivalences between the appropriately scaled filtered bases and $\mat T_\mathrm{s}$ and $\mat T_\mathrm{h}$. 
To numerically illustrate these equivalences, the spectra of these operators
and their preconditioned counterparts are illustrated in \Cref{fig:spectrum_Th_filtered_stars,fig:spectrum_Ts_filtered_loops}. These spectra correspond to a smoothly deformed sphere (see Fig.~\ref{fig:spectrum_Ts_filtered_loops} and \ref{fig:spectrum_Th_filtered_stars}), and the ordering of the singular values is obtained by projection against the graph Laplacians' eigenvectors. The original spectrum of $\mat T_\mathrm{s} $ and $\mat T_\mathrm{h}$ show the expected $\xi^{-1/2}$ and $\xi^{1/2}$---with $\xi$ the spectral index---behaviors, predicted by pseudo-differential operator theory. Given the construction of the preconditioners, it is then not surprising that the preconditioned operators show a spectrum bounded (and away from zero) with the expected variations in the spectrum.

\begin{figure}
    \centering
    \includegraphics{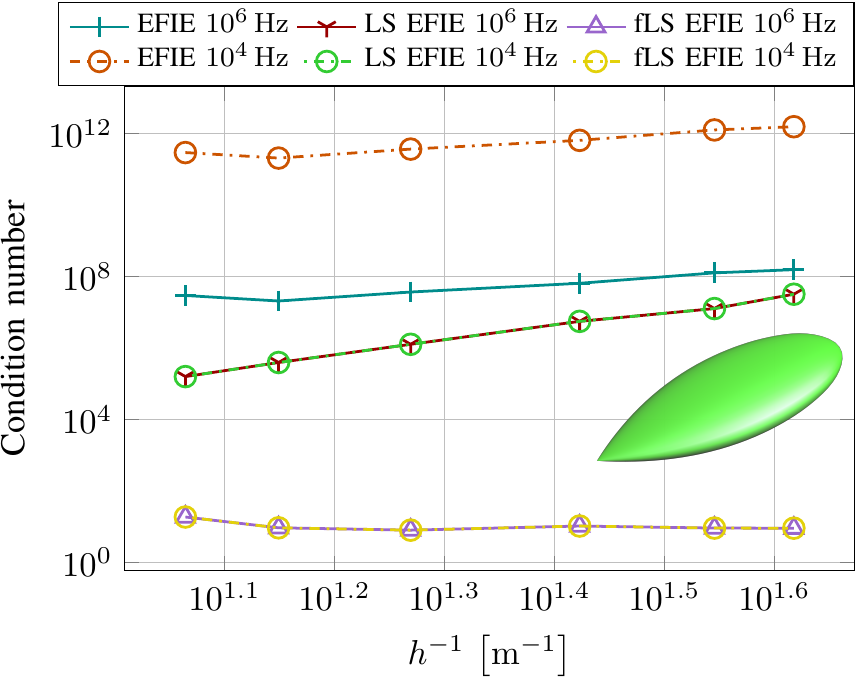}
    \caption{Condition number of the EFIE \eqref{eq:disc_EFIE}, Loop-Star EFIE, and filtered Loop-Star EFIE \eqref{eq:EFIE_wavelet} as a function of discretization for several frequencies. The condition number has been obtained after eliminating the isolated singular values, which have minimal impact on the convergence, arising from the deletion of one column from each of the preconditioning matrices.  The solid lines correspond to a simulating frequency of \SI{e6}{\hertz} and the dotted lines to a frequency of \SI{e4}{\hertz}. The simulated structure is the NASA almond re-scaled to be enclosed in a bounding box of diameter \SI{1.09}{\meter}. Perfect filters built out of SVD have been used in these results.}
    \label{fig:almond_LS_precond}
\end{figure}

To illustrate that the preconditioning schemes based on filtered bases do regularize the EFIE, the condition number of the original and preconditioned schemes will be compared for varying frequencies and discretizations.
First, the conditioning of a filtered Loop-Star preconditioned EFIE for the NASA almond \cite{wooEMProgrammerNotebookbenchmark1993} is reported in Figure~\ref{fig:almond_LS_precond}. The low frequency and dense discretization breakdowns of the original equations are apparent, while the preconditioned equation (corresponding to \eqref{eq:EFIE_wavelet}) shows a constant conditioning. This is in contrast with the standard Loop-Star approach that does regularize the low frequency conditioning breakdown, but actually worsens the dense discretization behavior of the equation.

\begin{figure}
    \centering
    \includegraphics{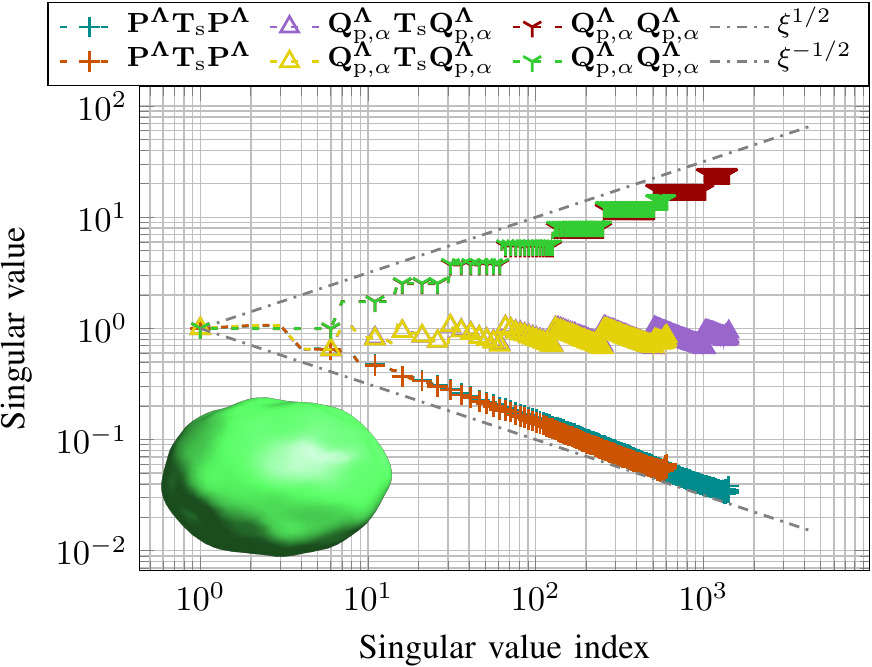}
    \caption{Spectrum of the solenoidal part of the vector potential, its preconditioner, and its preconditioned counterpart. These spectra have been obtained for a smoothly-deformed sphere with a maximum diameter of \SI{7.17}{\meter} (see insert), a frequency of \SI{e6}{\hertz}, and for two different average edge lengths \SI{0.31}{\meter} and \SI{0.20}{\meter}. The spectra have been normalized so that their first singular value is one, for readability. Perfect filters built out of SVD have been used in these results.}
    \label{fig:spectrum_Ts_filtered_proj_loops}
\end{figure}

\begin{figure}
    \centering
    \includegraphics{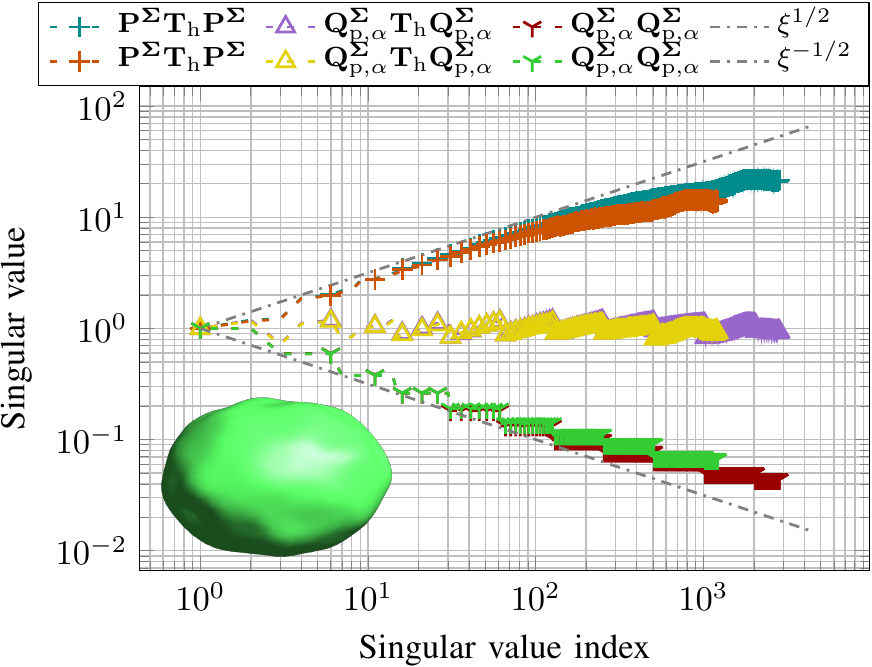}
    \caption{Spectrum of the non-solenoidal part of the scalar potential, its preconditioner, and its preconditioned counterpart. These spectra have been obtained for a smoothly-deformed sphere with a maximum diameter of \SI{7.17}{\meter} (see insert), a frequency of \SI{e6}{\hertz}, and for two different average edge lengths \SI{0.31}{\meter} and \SI{0.20}{\meter}. The spectra have been normalized so that their first singular value is one, for readability. Perfect filters built out of SVD have been used in these results.}
    \label{fig:spectrum_Th_filtered_proj_stars}
\end{figure}

\begin{figure}
    \centering
    \includegraphics{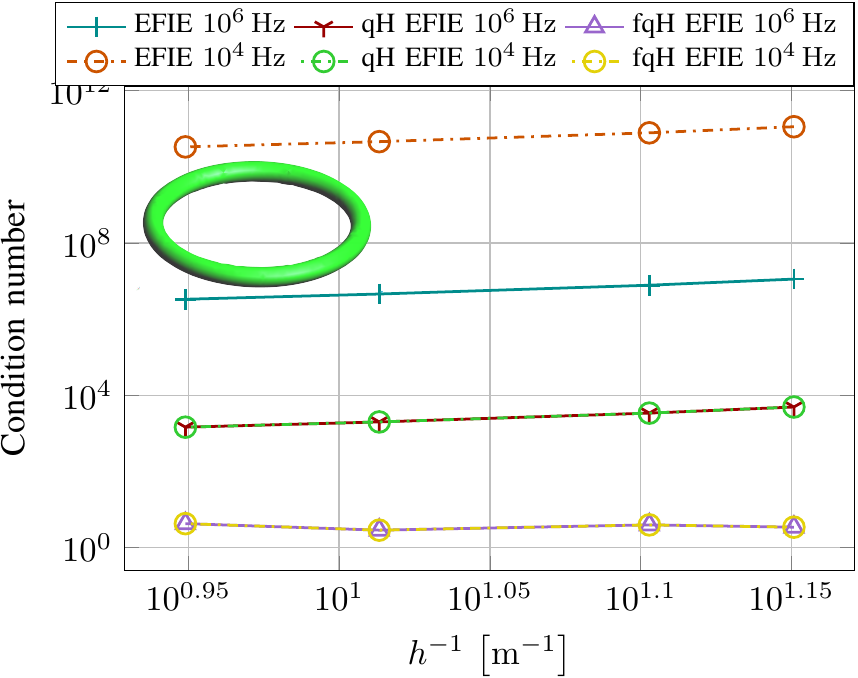}
    \caption{Condition number of the EFIE \eqref{eq:disc_EFIE}, quasi-Helmholtz (qH) projector EFIE, and filtered qH projector EFIE \eqref{eq:EFIE_proj} as a function of discretization for several frequencies. The solid lines correspond to a simulating frequency of \SI{e6}{\hertz} and the dotted lines to a frequency of \SI{e4}{\hertz}. The simulated structure is a torus with inner radius \SI{.9}{\meter} and outer radius \SI{1.1}{\meter}. Perfect filters built out of SVD have been used in these results.}
    \label{fig:torus_QH_precond_proj}
\end{figure}

A similar study has been performed with the filtered projectors schemes. In \Cref{fig:spectrum_Ts_filtered_proj_loops,fig:spectrum_Th_filtered_proj_stars} the spectra of the dominant solenoidal and non-solenoidal parts of the EFIE operators are displayed alongside their preconditioners.
%The stair-like behavior is also present in this graphs, as expected.
The preconditioning performance on the overall EFIE system is illustrated in \Cref{fig:torus_QH_precond_proj} for a torus.
The approach yields satisfactory conditioning that remains stable in both low frequency and dense discretization, which in turns shows that the scheme can also handle multiply-connected geometries. 


\begin{figure}
    \centering
    \includegraphics{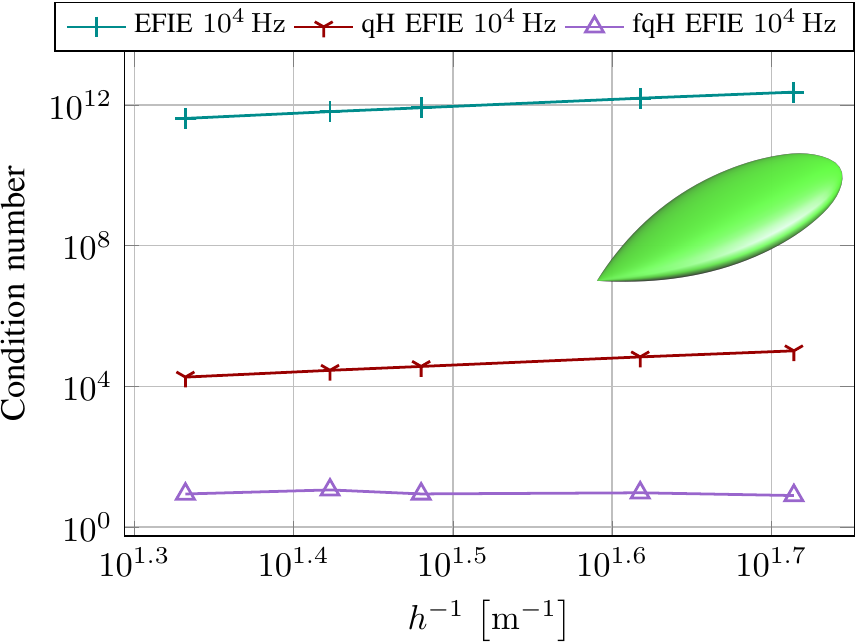}
    \caption{Condition number of the EFIE \eqref{eq:disc_EFIE}, quasi-Helmholtz (qH) projector EFIE, and filtered qH projector EFIE \eqref{eq:EFIE_proj} as a function of discretization for several frequencies. The simulated structure is the NASA almond re-scaled to be enclosed in a bounding box of diameter \SI{1.09}{\meter}. The preconditioner is built without using SVDs, but by leveraging Chebyshev-interpolated filters \eqref{eq:cheby_filters_def} corresponding to Butterworth filters of order \num{100}, expanded into \num{200} Chebyshev polynomials.}
    \label{fig:almond_QH_precond_proj_chebyshev}
\end{figure}

Finally, a conditioning study of the NASA almond is reported in \Cref{fig:almond_QH_precond_proj_chebyshev} that has been obtained using Chebyshev-interpolated filters \eqref{eq:cheby_filters_def} corresponding to Butterworth filters of order \num{100}, expanded into \num{200} Chebyshev polynomials. The coefficients of the filters are obtained via the norm estimates detailed in \eqref{eq:norm_star_proj_scaling} and \eqref{eq:norm_loop_proj_scaling} and the cutting point of the filters is determined using the approximate Laplacian spectrum described bellow \eqref{eq:factorized_butterworth}. The excellent stability of preconditioned scheme for a structure such as the NASA almond showcases the effectiveness of the scheme when using the fast techniques presented in this paper.

\section{Conclusion} \label{sec:conclusion}

A new family of strategies has been introduced for performing filtered quasi-Helmholtz decompositions of electromagnetic integral equations: the filtered Loop-Star decompositions and the quasi-Helmholtz Laplacian filters. These new tools are capable of manipulating large parts of the operators' spectra to obtain new families of preconditioners and fast solvers. A first application to the case of frequency and $h$-refinement preconditioning of the electric field integral equation has been presented and numerical results have shown the practical effectiveness of the newly proposed tools.

\appendices
\section{Complementarity of the projectors}
\label{app:complementarity}

In this appendix, we show that the properties $\Pn^{\Lambda} + \Pn^{\Sigma} = \mat I$ and $\Pnd^{\Lambda} + \Pnd^\Sigma =\mat I$ hold true on simply connected geometries. To this end, we first prove that the normalized coefficients $\vec{\tilde{j}}$ of the RWG functions can be decomposed with $\tmat \Lambda$ and $\tmat \Sigma$, similarly as in \eqref{eq:discreteHdecomposition} where we assume that the proper number of columns from the matrices have been removed as is standard to ensure a full column rank, such that
\begin{equation} \label{eq:ls_decomp_scaled}
    \vec{\tilde{j}} = \mat G^{-\frac{1}{2}} \vec j = \mat{\tilde{\Lambda}} \vec{\tilde{l}} + \mat{\tilde{\Sigma}} \vec{\tilde{s}}
\end{equation}
in which $\vec{\tilde{l}}$ and $\vec{\tilde{s}}$ are the coefficient vectors of the normalized Loop and Star parts in this decomposition. Since $\mat G$, $\mat{G_p}$, and $\mat{G}_\lambda$ are invertible matrices, we have $\text{rank}(\mat{\tilde{\Sigma}}) = \text{rank}(\mat{\Sigma})$ and $\text{rank}(\mat{\tilde{\Lambda}}) = \text{rank}(\mat{\Lambda})$. Moreover, since $\mat{\tilde{\Lambda}}^\T \mat{\tilde{\Sigma}} = \matO$, we also obtain that $\mat{\tilde{\Lambda}}$ and $\mat{\tilde{\Sigma}}$ have their column linearly independent, which yields $\text{rank}([\mat{\tilde{\Lambda}} ~ \mat{\tilde{\Sigma}}]) = \text{rank}(\mat{\tilde{\Lambda}} ) + \text{rank}(\mat{\tilde{\Sigma}} ) = N$, from which the existence and (unicity) of \eqref{eq:ls_decomp_scaled} follows.

Subsequently, using \eqref{eq:ls_decomp_scaled}, we can form a new set of normalized projectors to retrieve $\mat{\tilde{\Lambda}} \vec{\tilde{l}}$ and $\mat{\tilde{\Sigma}} \vec{\tilde{s}}$ separately. The first step is to apply $\mat{\tilde{\Lambda}}^\T$ and $\mat{\tilde{\Sigma}}^\T$ to \eqref{eq:ls_decomp_scaled} to express $\vec{\tilde{j}}$ in the two different bases 
\begin{gather}
    \mat{\tilde{\Lambda}}^\T \vec{\tilde{j}} = \mat{\tilde{\Lambda}}^\T  \mat{\tilde{\Lambda}} \vec{\tilde{l}}\,, \\
    \mat{\tilde{\Sigma}}^\T \vec{\tilde{j}} = \mat{\tilde{\Sigma}}^\T \mat{\tilde{\Sigma}} \vec{\tilde{s}}\,,
\end{gather}
since $\mat{\tilde{\Lambda}}^\T \mat{\tilde{\Sigma}} = \matO$ and $\mat{\tilde{\Sigma}}^\T \mat{\tilde{\Lambda}} = \matO$, given \eqref{eq:orthogonality1}. Subsequently, we express the coefficients of the normalized Loop and Star bases as a function of $\vec{\tilde{j}}$
\begin{gather}
    \vec{\tilde{l}} = \left( \mat{\tilde{\Lambda}}^\T  \mat{\tilde{\Lambda}} \right)^+ \mat{\tilde{\Lambda}}^\T \vec{\tilde{j}}\,, \label{eq:ltilde}\\
    \vec{\tilde{s}} = \left( \mat{\tilde{\Sigma}}^\T  \mat{\tilde{\Sigma}} \right)^+ \mat{\tilde{\Sigma}}^\T \vec{\tilde{j}}\,.
    \label{eq:stilde}
\end{gather}
Finally, we express $\mat{\tilde{\Lambda}}\vec{\tilde{l}}$ and $\mat{\tilde{\Sigma}} \vec{\tilde{s}}$ in terms of $\vec{\tilde{j}}$ by applying $\mat{\tilde{\Lambda}}$ and $\mat{\tilde{\Sigma}}$ to \eqref{eq:ltilde} and \eqref{eq:stilde}
\begin{gather}
    \mat{\tilde{\Lambda}} \vec{\tilde{l}} = \mat{\tilde{\Lambda}} \left(\mat{\tilde{\Lambda}}^\T  \mat{\tilde{\Lambda}} \right)^+ \mat{\tilde{\Lambda}}^\T \vec{\tilde{j}} = \Pn^{\Lambda} \vec{\tilde{j}}\,, \label{eq:Lltilde}\\
    \mat{\tilde{\Sigma}} \vec{\tilde{s}} = \mat{\tilde{\Sigma}} \left(\mat{\tilde{\Sigma}}^\T  \mat{\tilde{\Sigma}} \right)^+ \mat{\tilde{\Sigma}}^\T \vec{\tilde{j}} = \Pn^{\Sigma} \vec{\tilde{j}}\,, \label{eq:Sstilde}
\end{gather}
and we obtain that $\Pn^{\Lambda} + \Pn^{\Sigma} = \mat I$ by leveraging \eqref{eq:Lltilde}, \eqref{eq:Sstilde}, and \eqref{eq:ls_decomp_scaled}. Following the same procedure, except that now $\tilde{\dLambda}$ and $\tilde{\dSigma}$ are employed in the initial decomposition, we can show that the property $\Pnd^{\Lambda} + \Pnd^\Sigma =\mat I$ also holds true.

\section*{Acknowledgment}

This work has been funded in part by the European Research Council (ERC) under the European Union’s Horizon 2020 research and innovation program (ERC project 321, grant No.724846) and in part by the ANR Labex CominLabs under the project ``CYCLE''.

\ifCLASSOPTIONcaptionsoff
  \newpage
\fi

\bibliographystyle{IEEEtran}
\bibliography{MyLibrary,reference_adrien}


\begin{IEEEbiography}{Michael Shell}
Biography text here.
\end{IEEEbiography}

\begin{IEEEbiographynophoto}{John Doe}
Biography text here.
\end{IEEEbiographynophoto}

%\newpage

\begin{IEEEbiographynophoto}{Jane Doe}
Biography text here.
\end{IEEEbiographynophoto}


%\vfill


\end{document}

%% file: Packages.tex
\usepackage[utf8]{inputenc}

\usepackage[T1]{fontenc}
\usepackage{subcaption}
\usepackage{float}

\usepackage{cite}
\usepackage{url}
\usepackage{graphicx}
\usepackage{color}
\usepackage{placeins}
\usepackage{tabularx,colortbl}
\usepackage{ifthen}

\usepackage{pgf,tikz}
\usepackage{pgfplotstable}
\usepackage{pgfplots}
\usetikzlibrary{arrows,matrix,shadows,positioning,arrows.meta}
\usepackage{siunitx}
\pgfplotsset{compat=newest}
\usepackage{stmaryrd} % for double brackets
\usepackage{amsmath,amssymb,amsfonts}
\usepackage{amsthm}
\DeclareMathAlphabet{\mathbfsf}{\encodingdefault}{\sfdefault}{bx}{n}
\usepackage{bm}
\usepackage[OMLmathsfit]{isomath}

\usepackage{algorithmic}
\usepackage{textcomp}
\usepackage{xcolor}
\usepackage{mathtools}
\usepackage{todonotes}
\usepackage[belowskip=-15pt,aboveskip=0pt]{caption}
\usepackage{comment}
\usepackage[normalem]{ulem}
\usepackage{cleveref}

%% file: Macros.tex
\newcommand{\mr}{\mathrm}                   % Roman, non-italic
\newcommand{\vt}[1]{\boldsymbol{#1}} % Physical vectors
\newcommand{\vr}{\boldsymbol{r}} % Physical vectors
\newcommand{\veg}[1]{\bm{#1}}               % Bold for geometrical / physical vectors
\newcommand{\mat}[1]{\mathbf{#1}}       % NxN dimensional matrix
\newcommand{\matd}[1]{\mathbb{#1}}   % NxN dual dimensional matrix
\newcommand{\uv}[1]{\hat{\veg{#1}}}         % Unit vectors
\renewcommand{\vec}[1]{\mathsfbfit{#1}}     % N-dimensional vector
\newcommand{\op}[1]{\mathcal{#1}}           % Operator that takes a scalar-valued function
\newcommand{\vecop}[1]{\bm{\mathcal{#1}}}   % Operator that takes a vector-valued function
\newcommand{\vv}[1]{\vt{#1}}

\newcommand{\matI}{\mathbfsf{I}}
\newcommand{\matO}{\mathbfsf{0}}

\newcommand{\td}[1]{\todo[inline]{#1}}

\newcommand{\Pn}{ \tilde{\mat P}}
\newcommand{\Pnd}{\tilde{\matd P}}
\newcommand{\SigmaN}{\tilde{\mat \Sigma}}
\newcommand{\LambdaN}{\tilde{\mat \Lambda}}
\newcommand{\SigmaNd}{\tilde{\matd \Sigma}}
\newcommand{\LambdaNd}{\tilde{\matd \Lambda}}

\newcommand{\dd}{\mathrm{d}}  % The differential symbol is actually an operator

\newcommand{\im}{\mathrm{i}}  % For mathematician, Julia notation
\newcommand{\e}{\mathrm{e}}

\newcommand{\T}{\mr{T}}

\newcommand{\tmat}[1]{\mat{\tilde{#1}}}

%% file: preamble.tex
\pgfplotsset{compat=newest}
\pgfplotsset{plot coordinates/math parser=false}
\newlength\figureheight
\newlength\figurewidth
\pgfplotsset{every axis plot/.append style={line width=1.5pt},
    legend style={font=\footnotesize, 
        text height=1.0ex,
        draw=black,
        fill=white,
        legend cell align=left}}

\pgfplotsset{select coords between index/.style 2 args={
    x filter/.code={
        \ifnum\coordindex<#1\fi
        \ifnum\coordindex>#2\fi
    }
}}

\definecolor{amethyst}{rgb}{0.6, 0.4, 0.8}
\definecolor{mediumblue}{rgb}{0.0, 0.0, 0.8}
\definecolor{oucrimsonred}{rgb}{0.6, 0.0, 0.0}
\definecolor{darkcyan}{rgb}{0.0, 0.55, 0.55}

\definecolor{applegreen}{rgb}{0.55, 0.71, 0.0}
\definecolor{citrine}{rgb}{0.89, 0.82, 0.04}
\definecolor{limegreen}{rgb}{0.2, 0.8, 0.2}
\definecolor{burntorange}{rgb}{0.8, 0.33, 0.0}
\definecolor{gold(metallic)}{rgb}{0.83, 0.69, 0.22}

\pgfplotsset{
  compat=newest,
  GraphStyle/.style={
    width=\columnwidth, grid=both,legend style={font=\footnotesize}, mark
    size=0.3em, mark options={solid},
  },
  EFIEStyle/.style={
    thick, mark=+, mark options={solid}, dash pattern=on 2pt off 4pt, darkcyan
  },
  EFIEStyleComp/.style={
    thick, mark=+, mark options={solid}, dash pattern=on 2pt off 4pt, burntorange
  },
  LSEFIEStyle/.style={
    thick, mark=Mercedes star flipped, mark options={solid}, dash pattern=on 2pt off 4pt, dash phase=3pt, oucrimsonred
  },
  LSEFIEStyleComp/.style={
    thick, mark=Mercedes star flipped, mark options={solid}, dash pattern=on 2pt off 4pt, dash phase=3pt, limegreen
  },
  fLSEFIEStyle/.style={
    thick, mark=triangle, mark options={solid}, dash pattern=on 2pt off 4pt,dash phase=3pt, amethyst
  },
  fLSEFIEStyleComp/.style={
    thick, mark=triangle, mark options={solid}, dash pattern=on 2pt off 4pt,dash phase=3pt, citrine
  },
  PEFIEStyle/.style={
    thick, mark=diamond, mark options={solid}, dash pattern=on 2pt off 4pt, dash phase=3pt, mediumblue
  },
  PEFIEStyleComp/.style={
    thick, mark=diamond, mark options={solid}, dash pattern=on 2pt off 4pt, dash phase=3pt, gold(metallic)
  },
  TrendSqrtStyle/.style={
    thick, mark=none, dash dot, gray 
  },
  TrendInvSqrtStyle/.style={
    thick, mark=none, dash dot, black 
  },
}
\usepgfplotslibrary{units}

%% file: LapFiltersTAP.bbl
\begin{thebibliography}{10}
\providecommand{\url}[1]{#1}
\csname url@samestyle\endcsname
\providecommand{\newblock}{\relax}
\providecommand{\bibinfo}[2]{#2}
\providecommand{\BIBentrySTDinterwordspacing}{\spaceskip=0pt\relax}
\providecommand{\BIBentryALTinterwordstretchfactor}{4}
\providecommand{\BIBentryALTinterwordspacing}{\spaceskip=\fontdimen2\font plus
\BIBentryALTinterwordstretchfactor\fontdimen3\font minus
  \fontdimen4\font\relax}
\providecommand{\BIBforeignlanguage}[2]{{%
\expandafter\ifx\csname l@#1\endcsname\relax
\typeout{** WARNING: IEEEtran.bst: No hyphenation pattern has been}%
\typeout{** loaded for the language `#1'. Using the pattern for}%
\typeout{** the default language instead.}%
\else
\language=\csname l@#1\endcsname
\fi
#2}}
\providecommand{\BIBdecl}{\relax}
\BIBdecl

\bibitem{nedelec_acoustic_2001}
J.-C. Nedelec, \emph{Acoustic and {{Electromagnetic Equations}}: {{Integral
  Representations}} for {{Harmonic Problems}}}.\hskip 1em plus 0.5em minus
  0.4em\relax {Springer Science \& Business Media}, Mar. 2001.

\bibitem{gibson_method_2014}
W.~C. Gibson, \emph{The {{Method}} of {{Moments}} in {{Electromagnetics}} 2nd
  {{Ed}}.}\hskip 1em plus 0.5em minus 0.4em\relax {CRC press}, 2014.

\bibitem{jin_theory_2015}
J.-M. Jin, \emph{Theory and Computation of Electromagnetic Fields},
  2nd~ed.\hskip 1em plus 0.5em minus 0.4em\relax {Piscataway, NJ}: {IEEE
  Press}, 2015.

\bibitem{chew_fast_2001}
W.~C. Chew, J.-M. Jin, E.~Michielssen, and J.~M. Song, Eds., \emph{Fast and
  {{Efficient Algorithms}} in {{Computational Electromagnetics}}}.\hskip 1em
  plus 0.5em minus 0.4em\relax {Artech House}, 2001.

\bibitem{axelsson_iterative_1996}
O.~Axelsson, \emph{Iterative Solution Methods}.\hskip 1em plus 0.5em minus
  0.4em\relax {Cambridge university press}, 1996.

\bibitem{colton_integral_2013}
D.~L. Colton and R.~Kress, \emph{Integral Equation Methods in Scattering
  Theory}.\hskip 1em plus 0.5em minus 0.4em\relax {Philadelphia}: {Society for
  Industrial and Applied Mathematics}, 2013.

\bibitem{adrianElectromagneticIntegralEquations2021}
S.~B. Adrian, A.~D{\'e}ly, D.~Consoli, A.~Merlini, and F.~P. Andriulli,
  ``Electromagnetic {{Integral Equations}}: {{Insights}} in {{Conditioning}}
  and {{Preconditioning}},'' \emph{IEEE Open Journal of Antennas and
  Propagation}, vol.~2, pp. 1143--1174, 2021.

\bibitem{wilton_topological_1983}
D.~R. Wilton, ``Topological consideration in surface patch and volume cell
  modeling of electromagnetic scatterers,'' in \emph{Proc. {{URSI Int}}.
  {{Symp}}. {{Electromagn}}. {{Theory}}}, 1983, pp. 65--68.

\bibitem{vecchi_loop-star_1999}
G.~Vecchi, ``Loop-star decomposition of basis functions in the discretization
  of the {{EFIE}},'' \emph{IEEE Transactions on Antennas and Propagation},
  vol.~47, no.~2, pp. 339--346, 1999.

\bibitem{lee_loop_2003}
J.-F. Lee, R.~Lee, and R.~Burkholder, ``Loop star basis functions and a robust
  preconditioner for {{EFIE}} scattering problems,'' \emph{IEEE Transactions on
  Antennas and Propagation}, vol.~51, no.~8, pp. 1855--1863, Aug. 2003.

\bibitem{eibert_iterative-solver_2004}
T.~F. Eibert, ``Iterative-solver convergence for loop-star and loop-tree
  decompositions in method-of-moments solutions of the electric-field integral
  equation,'' \emph{Antennas and Propagation Magazine, IEEE}, vol.~46, no.~3,
  pp. 80--85, 2004.

\bibitem{andriulli_loop-star_2012}
F.~P. Andriulli, ``Loop-{{Star}} and {{Loop-Tree Decompositions}}: {{Analysis}}
  and {{Efficient Algorithms}},'' \emph{IEEE Transactions on Antennas and
  Propagation}, vol.~60, no.~5, pp. 2347--2356, May 2012.

\bibitem{vipianaMultiresolutionMethodMoments2005}
F.~Vipiana, P.~Pirinoli, and G.~Vecchi, ``A multiresolution method of moments
  for triangular meshes,'' \emph{IEEE transactions on antennas and
  propagation}, vol.~53, no.~7, pp. 2247--2258, 2005.

\bibitem{andriulli_multiresolution_2007}
F.~P. Andriulli, A.~Tabacco, and G.~Vecchi, ``A multiresolution approach to the
  electric field integral equation in antenna problems,'' \emph{SIAM Journal on
  Scientific Computing}, vol.~29, no.~1, pp. 1--21, 2007.

\bibitem{andriulli_hierarchical_2008}
F.~P. Andriulli, F.~Vipiana, and G.~Vecchi, ``Hierarchical {{Bases}} for
  {{Nonhierarchic}} 3-{{D Triangular Meshes}},'' \emph{IEEE Transactions on
  Antennas and Propagation}, vol.~56, no.~8, pp. 2288--2297, Aug. 2008.

\bibitem{chenMultiresolutionCurvilinearRao2009}
R.-S. Chen, J.~Ding, D.~Ding, Z.~Fan, and D.~Wang, ``A multiresolution
  curvilinear rao--wilton--glisson basis function for fast analysis of
  electromagnetic scattering,'' \emph{IEEE Transactions on Antennas and
  Propagation}, vol.~57, no.~10, pp. 3179--3188, 2009.

\bibitem{guzman_hierarchical_2016-1}
J.~E.~O. Guzman, S.~Adrian, R.~Mitharwal, Y.~Beghein, T.~Eibert, K.~Cools, and
  F.~Andriulli, ``On the {{Hierarchical Preconditioning}} of the {{PMCHWT
  Integral Equation}} on {{Simply}} and {{Multiply Connected Geometries}},''
  \emph{IEEE Antennas and Wireless Propagation Letters}, vol.~PP, no.~99, pp.
  1--1, 2016.

\bibitem{adrian_hierarchical_2017}
S.~Adrian, F.~Andriulli, and T.~Eibert, ``A hierarchical preconditioner for the
  electric field integral equation on unstructured meshes based on primal and
  dual {{Haar}} bases,'' \emph{Journal of Computational Physics}, vol. 330, pp.
  365--379, Feb. 2017.

\bibitem{adams_stabilisation_1999}
R.~Adams and G.~Brown, ``Stabilisation procedure for electric field integral
  equation,'' \emph{Electronics Letters}, vol.~35, no.~23, pp. 2015--2016, Nov.
  1999.

\bibitem{contopanagos_well-conditioned_2002}
H.~Contopanagos, B.~Dembart, M.~Epton, J.~Ottusch, V.~Rokhlin, J.~Visher, and
  S.~Wandzura, ``Well-conditioned boundary integral equations for
  three-dimensional electromagnetic scattering,'' \emph{IEEE Transactions on
  Antennas and Propagation}, vol.~50, no.~12, pp. 1824--1830, Dec. 2002.

\bibitem{adams_physical_2004}
R.~Adams, ``Physical and {{Analytical Properties}} of a {{Stabilized Electric
  Field Integral Equation}},'' \emph{IEEE Transactions on Antennas and
  Propagation}, vol.~52, no.~2, pp. 362--372, Feb. 2004.

\bibitem{adams_numerical_2004}
R.~Adams and N.~Champagne, ``A {{Numerical Implementation}} of a {{Modified
  Form}} of the {{Electric Field Integral Equation}},'' \emph{IEEE Transactions
  on Antennas and Propagation}, vol.~52, no.~9, pp. 2262--2266, Sep. 2004.

\bibitem{borel_new_2005}
S.~Borel, D.~Levadoux, and F.~Alouges, ``A new well-conditioned {{Integral}}
  formulation for {{Maxwell}} equations in three dimensions,'' \emph{IEEE
  Transactions on Antennas and Propagation}, vol.~53, no.~9, pp. 2995--3004,
  Sep. 2005.

\bibitem{andriulli_multiplicative_2008}
F.~P. Andriulli, K.~Cools, H.~Bagci, F.~Olyslager, A.~Buffa, S.~Christiansen,
  and E.~Michielssen, ``A {{Multiplicative Calderon Preconditioner}} for the
  {{Electric Field Integral Equation}},'' \emph{IEEE Transactions on Antennas
  and Propagation}, vol.~56, no.~8, pp. 2398--2412, Aug. 2008.

\bibitem{stephanson_preconditioned_2009}
M.~B. Stephanson and J.-F. Lee, ``Preconditioned {{Electric Field Integral
  Equation Using Calderon Identities}} and {{Dual Loop}}/{{Star Basis
  Functions}},'' \emph{IEEE Transactions on Antennas and Propagation}, vol.~57,
  no.~4, pp. 1274--1279, Apr. 2009.

\bibitem{borel2005new}
S.~Borel, D.~P. Levadoux, and F.~Alouges, ``A new well-conditioned integral
  formulation for maxwell equations in three dimensions,'' \emph{IEEE
  transactions on antennas and propagation}, vol.~53, no.~9, pp. 2995--3004,
  2005.

\bibitem{andriulli_well-conditioned_2013}
F.~P. Andriulli, K.~Cools, I.~Bogaert, and E.~Michielssen, ``On a
  {{Well-Conditioned Electric Field Integral Operator}} for {{Multiply
  Connected Geometries}},'' \emph{IEEE Transactions on Antennas and
  Propagation}, vol.~61, no.~4, pp. 2077--2087, Apr. 2013.

\bibitem{dobbelaere_calderon_2015}
D.~Dobbelaere, D.~De~Zutter, J.~Van~Hese, J.~Sercu, T.~Boonen, and H.~Rogier,
  ``A {{Calder\'on}} multiplicative preconditioner for the electromagnetic
  {{Poincar\'e}}\textendash{{Steklov}} operator of a heterogeneous domain with
  scattering applications,'' \emph{Journal of Computational Physics}, vol. 303,
  pp. 355--371, Dec. 2015.

\bibitem{merlini_magnetic_2020}
A.~Merlini, Y.~Beghein, K.~Cools, E.~Michielssen, and F.~P. Andriulli,
  ``Magnetic and {{Combined Field Integral Equations Based}} on the
  {{Quasi-Helmholtz Projectors}},'' \emph{IEEE Transactions on Antennas and
  Propagation}, vol.~68, no.~5, pp. 3834--3846, May 2020.

\bibitem{adrian_refinement-free_2018}
S.~Adrian, F.~Andriulli, and T.~Eibert, ``On a {{Refinement-Free Calder\'on
  Multiplicative Preconditioner}} for the {{Electric Field Integral
  Equation}},'' \emph{Journal of Computational Physics}, Oct. 2018.

\bibitem{johnsonPolynomialPreconditionersConjugate1983}
O.~G. Johnson, C.~A. Micchelli, and G.~Paul, ``Polynomial preconditioners for
  conjugate gradient calculations,'' \emph{SIAM Journal on Numerical Analysis},
  vol.~20, no.~2, pp. 362--376, 1983.

\bibitem{ashbyComparisonAdaptiveChebyshev1992}
S.~F. Ashby, T.~A. Manteuffel, and J.~S. Otto, ``A comparison of adaptive
  chebyshev and least squares polynomial preconditioning for hermitian positive
  definite linear systems,'' \emph{SIAM Journal on Scientific and Statistical
  Computing}, vol.~13, no.~1, pp. 1--29, 1992.

\bibitem{hammond_wavelets_2011}
D.~K. Hammond, P.~Vandergheynst, and R.~Gribonval, ``Wavelets on graphs via
  spectral graph theory,'' \emph{Applied and Computational Harmonic Analysis},
  vol.~30, no.~2, pp. 129--150, Mar. 2011.

\bibitem{levie2018cayleynets}
R.~Levie, F.~Monti, X.~Bresson, and M.~M. Bronstein, ``Cayleynets: Graph
  convolutional neural networks with complex rational spectral filters,''
  \emph{IEEE Transactions on Signal Processing}, vol.~67, no.~1, pp. 97--109,
  2018.

\bibitem{rahmouni_new_2019}
L.~Rahmouni and F.~P. Andriulli, ``A {{New Preconditioner}} for the {{EFIE
  Based}} on {{Primal}} and {{Dual Graph Laplacian Spectral Filters}},'' in
  \emph{2019 {{International Conference}} on {{Electromagnetics}} in {{Advanced
  Applications}} ({{ICEAA}})}.\hskip 1em plus 0.5em minus 0.4em\relax {Granada,
  Spain}: {IEEE}, Sep. 2019, pp. 1342--1344.

\bibitem{9886850}
A.~Merlini, C.~Henry, D.~Consoli, L.~Rahmouni, and F.~P. Andriulli, ``Laplacian
  filters for integral equations: Further developments and fast algorithms,''
  in \emph{2022 IEEE International Symposium on Antennas and Propagation and
  USNC-URSI Radio Science Meeting (AP-S/URSI)}, 2022, pp. 1932--1933.

\bibitem{rao_electromagnetic_1982}
S.~Rao, D.~Wilton, and A.~Glisson, ``Electromagnetic scattering by surfaces of
  arbitrary shape,'' \emph{IEEE Transactions on Antennas and Propagation},
  vol.~30, no.~3, pp. 409--418, May 1982.

\bibitem{buffa_dual_2007}
A.~Buffa and S.~Christiansen, ``A dual finite element complex on the
  barycentric refinement,'' \emph{Mathematics of Computation}, vol.~76, no.
  260, pp. 1743--1769, 2007.

\bibitem{chen_electromagnetic_1990}
Q.~Chen and D.~Wilton, ``Electromagnetic scattering by three-dimensional
  arbitrary complex material/conducting bodies,'' in \emph{Antennas and
  {{Propagation Society International Symposium}}, 1990. {{AP-S}}. {{Merging
  Technologies}} for the 90's. {{Digest}}.}, May 1990, pp. 590--593 vol.2.

\bibitem{mautzEfieldSolutionConducting1984}
J.~Mautz and R.~Harrington, ``An {{E-field}} solution for a conducting surface
  small or comparable to the wavelength,'' \emph{IEEE Transactions on Antennas
  and Propagation}, vol.~32, no.~4, pp. 330--339, 1984.

\bibitem{limNovelTechniqueCalculate}
J.~S. Lim, S.~Rao, and D.~R. Wilton, ``A {{Novel Technique}} to {{Calculate}}
  the {{Electromagnetic Scattering}} by {{Surfaces}} of {{Arbitrary Shape}},''
  in \emph{1993 {{URSI Radio Science Meeting Digest}}}, 1993, p. 322.

\bibitem{wuStudyTwoNumerical1995}
W.-L. Wu, A.~W. Glisson, and D.~Kajfez, ``A study of two numerical solution
  procedures for the electric field integral equation at low frequency,''
  \emph{Applied Computational Electromagnetics Society Journal}, vol.~10,
  no.~3, pp. 69--80, 1995.

\bibitem{golub_matrix_2012}
G.~H. Golub and C.~F. Van~Loan, \emph{Matrix Computations}.\hskip 1em plus
  0.5em minus 0.4em\relax {JHU Press}, 2012, vol.~3.

\bibitem{higham_functions_2008}
N.~Higham, \emph{Functions of {{Matrices}}}, ser. Other {{Titles}} in {{Applied
  Mathematics}}.\hskip 1em plus 0.5em minus 0.4em\relax {Society for Industrial
  and Applied Mathematics}, Jan. 2008.

\bibitem{press_numerical_2007}
W.~H. Press, \emph{Numerical {{Recipes}} 3rd {{Edition}}: {{The Art}} of
  {{Scientific Computing}}}.\hskip 1em plus 0.5em minus 0.4em\relax {Cambridge
  University Press}, Sep. 2007.

\bibitem{mitharwal_multiplicative_2014}
R.~Mitharwal and F.~Andriulli, ``On the {{Multiplicative Regularization}} of
  {{Graph Laplacians}} on {{Closed}} and {{Open Structures With Applications}}
  to {{Spectral Partitioning}},'' \emph{IEEE Access}, vol.~2, pp. 788--796,
  2014.

\bibitem{oneilSecondkindIntegralEquations2018}
M.~O'Neil, ``Second-kind integral equations for the {{Laplace-Beltrami}}
  problem on surfaces in three dimensions,'' \emph{Advances in Computational
  Mathematics}, vol.~44, no.~5, pp. 1385--1409, Oct. 2018.

\bibitem{zhao_integral_2000}
J.-S. Zhao and W.~C. Chew, ``Integral equation solution of {{Maxwell}}'s
  equations from zero frequency to microwave frequencies,'' \emph{Antennas and
  Propagation, IEEE Transactions on}, vol.~48, no.~10, pp. 1635--1645, 2000.

\bibitem{arnold2007compatible}
D.~N. Arnold, P.~B. Bochev, R.~B. Lehoucq, R.~A. Nicolaides, and M.~Shashkov,
  \emph{Compatible spatial discretizations}.\hskip 1em plus 0.5em minus
  0.4em\relax Springer Science \& Business Media, 2007, vol. 142.

\bibitem{boubendir_well-conditioned_2014}
Y.~Boubendir and C.~Turc, ``Well-conditioned boundary integral equation
  formulations for the solution of high-frequency electromagnetic scattering
  problems,'' \emph{Computers \& Mathematics with Applications}, vol.~67,
  no.~10, pp. 1772--1805, Jun. 2014.

\bibitem{wooEMProgrammerNotebookbenchmark1993}
A.~Woo, H.~Wang, M.~Schuh, and M.~Sanders, ``{{EM}} programmer's
  notebook-benchmark radar targets for the validation of computational
  electromagnetics programs,'' \emph{IEEE Antennas and Propagation Magazine},
  vol.~35, no.~1, pp. 84--89, 1993.

\end{thebibliography}
